\documentclass[3p,times]{elsarticle}

\usepackage{my_ecrc}

\runauth{R. P. Agarwal, D. Baleanu, J. J. Nieto, D. F. M. Torres and Y. Zhou}

\usepackage{amssymb}
\usepackage{amsmath,amsthm}
\usepackage{pifont,color,subfigure}
\usepackage[figuresright]{rotating}

\usepackage{latexsym}
\usepackage{mathrsfs}

\usepackage{natbib}
\usepackage[colorlinks=true,linkcolor=black, citecolor=blue, urlcolor=blue]{hyperref}
	
\newtheorem{theorem}{Theorem}[section]
\newtheorem{corollary}{Corollary}[section]

\newtheorem{definition}{Definition}[section]
\newtheorem{remark}{Remark}[section]

\theoremstyle{definition} \theoremstyle{remark}

\begin{document}

\begin{frontmatter}

\dochead{\small [This is a preprint of a paper whose final and definite form is with 
\emph{Journal of Computational and Applied Mathematics}, ISSN:~0377-0427. 
Submitted 17-July-2017; Revised 18-Sept-2017; Accepted for publication 20-Sept-2017.]}

\title{A survey on fuzzy fractional differential and optimal control nonlocal evolution equations}

\author[add1]{Ravi P. Agarwal}
\ead{ravi.agarwal@tamuk.edu}

\author[add2,add3]{Dumitru Baleanu}
\ead{dumitru@cankaya.edu.tr}

\author[add4]{Juan J. Nieto}
\ead{juanjose.nieto.roig@usc.es}

\author[add5]{Delfim F. M. Torres\corref{cor1}}
\cortext[cor1]{Corresponding author}
\ead{delfim@ua.pt}

\author[add6,add7]{Yong Zhou}
\ead{yzhou@xtu.edu.cn}


\address[add1]{Texas A \& M University--Kingvsille, 
Department of Mathematics, Kingsville, TX 78363}

\address[add2]{Department of Mathematics,
Cankaya University, Ankara, Turkey}

\address[add3]{Institute of Space Sciences, Magurele--Bucharest, Romania}

\address[add4]{Facultade de Matem\'{a}ticas, Universidade de Santiago de Compostela, 
15782 Santiago de Compostela, Spain}

\address[add5]{Center for Research and Development in Mathematics and Applications (CIDMA),\\
Department of Mathematics, University of Aveiro, 3810-193 Aveiro, Portugal}

\address[add6]{Faculty of Mathematics and Computational Science, Xiangtan University,
Hunan 411105, P. R. China}

\address[add7]{Nonlinear Analysis and Applied Mathematics (NAAM) Research Group,\\ 
Faculty of Science, King Abdulaziz University, Jeddah 21589, Saudi Arabia}


\begin{abstract}
We survey some representative results on fuzzy fractional differential equations,
controllability, approximate controllability, optimal control, and optimal feedback control 
for several different kinds of fractional evolution equations.
Optimality and relaxation of multiple control problems,
described by nonlinear fractional differential equations
with nonlocal control conditions in Banach spaces, are considered.
\end{abstract}




\end{frontmatter}


\section{Introduction}

Memory and hereditary properties of different materials 
and processes in electrical circuits, biology, biomechanics,
electrochemistry, control, porous media and electromagnetic processes, 
are widely recognized to be well predicted by using fractional 
differential operators \cite{AMA.4,AMA.13,32,29,AMA.39}.
During the last decades, the subject of fractional calculus, 
and its potential applications, have gained an increase of importance, 
mainly because it has become a powerful tool with accurate and successful results
in modeling several complex phenomena in numerous seemingly
diverse and widespread fields of science and engineering
\cite{1,AMA.16,AMA.34,AMA.38}. Fractional calculus is 
not only a productive and emerging field, it also represents 
a new philosophy how to construct and apply a certain type 
of nonlocal operators to real world problems. The ones possessing 
both nonlocal effects as well as uncertainty behaviors represent 
an interesting phenomena.  Researchers started to combine, in an 
intelligent way, the notions of fractional with fuzzy, therefore 
a hybrid operator called fuzzy fractional operator emerges. 
In this manuscript, we begin by presenting a review on fractional 
differential equations under uncertainty.

In one of the earliest works, Agarwal et al. \cite{1} took the initiative 
and introduced fuzzy fractional calculus to handle fractional-order systems 
with uncertain initial values or uncertain relationships between parameters. 
Arshad and Lupulescu \cite{2} utilized the results reported in \cite{1} 
and they proved the existence and uniqueness of fractional differential 
equations with uncertainty. Afterward, Allahviranloo et al. \cite{3} 
employed the Riemann--Liouville generalized H-differentiability in order 
to solve the fuzzy fractional differential equations (FFDEs) and presented 
some new results under this notion. Salahshour et al. \cite{4} 
apply the technique of fuzzy Laplace transforms and solved some types 
of FFDEs based on the Riemann--Liouville fuzzy derivative.
Based on the delta-Hukuhara derivative for fuzzy valued functions,
Fard et al. established stability criteria for hybrid fuzzy
systems on time scales in the Lyapunov sense \cite{MyID:341}.
In \cite{MyID:363}, Fard et al. solve a class
of fuzzy fractional optimal control problems, 
where the coefficients of the system can be time-dependent. M
ore precisely, they establish a weak version of the Pontryagin maximum
principle for fuzzy fractional optimal control problems 
depending on generalized Hukuhara fractional Caputo derivatives \cite{MyID:363}.

Generally, the majority of the FFDEs as same as FDEs do not have exact solutions. 
As a result, approximate and numerical procedures are important to be developed
\cite{MyID:264}. On the other hand, because many of the parameters in 
mathematical models often do not appear explicitly,
modeling of natural phenomena using fuzzy fractional models 
plays an important role in various disciplines. Hence, 
it motivates the researchers to investigate effective numerical methods 
with error analysis to approximate the FFDEs.  As a result, researchers 
started to develop numerical techniques for FFDEs.  
Mazandarani and Vahidian Kamyad \cite{5} introduced a fuzzy approximate 
solution using the Euler method to solve FFDEs. Ahmadian et al. \cite{6} 
adopted the operational Jacobi operational matrix based on the fuzzy 
Caputo fractional derivative using shifted Jacobi polynomials. 
The clear advantage of the usage of this method is that the matrix 
operators have the main role to find the approximate fuzzy solution 
of FFDEs instead of considering the methods required the complicated 
fractional derivatives and their calculations.

Ghaemi et al. \cite{7} adapted a spectral method for the numerical 
solution of fuzzy fractional kinetic equations. The proposed method 
is characterized by its simplicity, efficiency, and high accuracy. 
Using the proposed method, they could reach a suitable approximation 
of the amount of the concentration value of xylose after a determined time 
that is important to analyze the kinetic data in the chemical process.  
Ahmadian et al. \cite{8} exploited a cluster of orthogonal functions, 
named shifted Legendre functions, to solve FFDEs under Caputo type. 
The benefit of the shifted Legendre operational matrices method, 
over other existing orthogonal polynomials, is its simplicity 
of execution as well as some other advantages. The achieved solutions 
present satisfactory results, obtained with only a small number 
of Legendre polynomials.

Fuzzy theory provides a suitable way to objectively account 
for parameter uncertainty in models. Fuzzy logic approaches 
appear promising in preclinical applications and might be useful 
in drug discovery and design. In this regards, Ahmadian et al. \cite{9} 
developed a tau method based on the Jacobi operational matrix 
to numerically solve the pharmacokinetics-pharmacodynamic equation, 
arising from drug assimilation into the bloodstream. The comparison 
of the results shows that the present method is a powerful mathematical 
tool for finding the numerical solutions of a generalized linear 
fuzzy fractional pharmacokinetics-pharmacodynamic equation. 
Balooch Shahriyar et al. \cite{10} investigated an analytical method 
(eigenvalue-eigenvector) for solving a system of FFDEs under fuzzy 
Caputo's derivative. To this end, they exploited generalized 
H-differentiability and derived the solutions based on this concept.
Ahmadian et al. \cite{11}  were confined with the application of 
Legendre operational matrix for solving FFDEs arising in the drug 
delivery model into the bloodstream. The main motivation of this research 
is to recommend a suitable way to approximate fuzzy fractional 
pharmacokinetics-pharmacodynamic models using a shifted Legendre tau approach. 
This strategy demands a formula for fuzzy fractional-order Caputo derivatives 
of shifted Legendre polynomials.

Mazandarani and Najariyan \cite{12} introduced two definitions 
of differentiability of type-2 fuzzy number valued functions 
of fractional order. The definitions are in the sense of Riemann--Liouville 
and Caputo. The methods, under type-2 fuzzy sets theory, will lead to an 
increase in the computational cost, although it is closer 
to the originality of the model. Salahshour et al. \cite{13} developed 
the notion of Caputo's H-differentiability, based on the generalized 
Hukuhara difference, to solve the FFDE.  To this end, they revisited 
Caputo's derivatives, and proposed novel fuzzy Laplace transforms 
and their inverses, with an analytical method to tackle the deficiencies 
in the state-of-the-art methods. Experimental results using some 
real-world problems (nuclear decay equation and Basset problem, 
illustrated the effectiveness and applicability of the proposed method). 
Simultaneously, the authors in \cite{121} investigated 
an effective numerical method with error analysis to approximate the fuzzy
time-fractional Bloch equations (FTFBE) \eqref{eq21} on the time interval
$J = (0, T]$, with a view to be employed in the image processing
domain in near time. 

Employing Laplace transforms, the authors in \cite{36} proposed 
a novel efficient technique for the solution of FFDEs that can 
efficiently make the original problem easier to achieve the numerical solution.
The suggested algorithm for the FFDEs use the fuzzy fractional
derivative of Caputo type in the range of $\alpha\in (0, 1]$ and is
potentially useful in solving fractional viscoelastic problems
under uncertainty. Chehlabi and Allahviranloo \cite{37} studied 
fuzzy linear fractional differential equations of order $0 < \alpha \leq 1$ 
under Riemann--Liouville H-differentiability. Also, it is corrected  
some previous results and obtained new solutions by using 
fractional hyperbolic functions and their properties.

There has been a significant development in nonlocal problems for fractional
differential equations or inclusions: see, for instance,
\cite{AMA.6,AMA.8,AMA.9,MyID:264,AMA.33,MR2824730,AMA.46}.
Indeed, nonlinear fractional differential equations have, in recent years,
been object of an increasing interest because of their wide applicability
in nonlinear oscillations of earthquakes, many physical phenomena such
as seepage flow in porous media, and in fluid dynamic traffic model
\cite{AMA.23,AMA.25,AMA.26}. On the other hand,
there could be no manufacturing, no vehicles,
no computers, and no regulated environment, without control systems.
Control systems are most often based on the principle of feedback,
whereby the signal to be controlled is compared to a desired reference
signal and the discrepancy used to compute corrective control actions \cite{AMA.18}.

The idea of controllability is an essential characteristic to a control 
framework exhibiting numerous control issues such as adjustment of unsteady 
frameworks by input control. Recently, control issues have been addressed  
by many physicists, engineers and mathematicians, and significant contribution 
on theoretical and application aspects of the topic can be found 
in the related literature \cite{MR3316531}.

As is well-known, the problems of exact and approximate
controllability are to be distinguished \cite{MR3244467}. 
In general, in infinite dimensional spaces, 
the concept of exact controllability is usually too strong. 
Therefore, the class of evolution equations consisting 
of fractional diffusion equations must be treated by the weaker concept of
controllability, namely approximate controllability \cite{MR3648947}. 
Recently, many works pay attention to study approximate controllability 
of different types of fractional evolution systems \cite{MyID:264,MR3244467,MR3291089}.

Over the last years, one of the fields of science that has been well established
is the fractional calculus of variations: see \cite{AMA.30,AMA.29,AMA.24} and references therein.
Moreover, a generalization of this area, namely the fractional optimal control,
is a topic of research by many authors \cite{AMA.2,AMA.19}.
The fractional optimal control of a distributed system is an optimal control problem
for which the system dynamics is defined with partial fractional differential 
equations \cite{AMA.35,MR3679433}. The calculus of variations, with constraints being 
sets of solutions of control systems, allow us to justify, while performing 
numerical calculations, the passage from a nonconvex optimal control problem 
to the convexified optimal control problem. We then approximate the latter problem 
by a sequence of smooth and convex optimal control
problems, for which the optimality conditions are known
and methods of their numerical resolution are well developed.

The delay evolution systems is an important class of distributed parameter systems, 
and optimal control of infinite dimensional systems is a remarkable subject 
in control theory \cite{MR3260698,MR2342658,MR2980253}.
In the last years, fractional evolution systems in
infinite dimensional spaces attracted many authors.
When the fractional differential equations describe the
performance index and system dynamics, an optimal control problem
reduces to a fractional optimal control problem \cite{MR3343434,MR3124694}. 
The fractional optimal control of a distributed system is a fractional optimal
control for which system dynamics are defined with partial
fractional differential equations. There has been very little work
in the area of fractional optimal control problem in infinite
dimensional spaces, especially optimal controls of fractional finite time delay
evolution system. See Sections~\ref{sec:exist_oc} and \ref{sec7}.

Sobolev type semilinear equations serve as an abstract formulation
of partial differential equations, which arise in various applications such
as in the flow of fluid through fissured rocks, thermodynamics, and shear
in second order fluids. Further, the fractional differential equations
of Sobolev type appear in the theory of control of dynamical systems,
when the controlled system and/or the controller is described
by a fractional differential equation of Sobolev type. Furthermore,
the mathematical modeling and simulations of systems and processes
are based on the description of their properties in terms of fractional
differential equations of Sobolev type. These new models are more adequate
than previously used integer order models, so fractional order differential
equations of Sobolev type have been investigated by many researchers:
see, for example, Fe\u ckan, Wang and Zhou \cite{AMA.17}
and Li, Liang and Xu \cite{AMA.27}.
In \cite{MR3260698,MR3244467}, the notion of nonlocal control condition 
is introduced and  a new kind of Sobolev type condition presented.
Kamocki \cite{AMA.22} studied the existence of optimal solutions 
to fractional optimal control problems. Liu et al. \cite{AMA.28} 
established the relaxation for nonconvex optimal control problems 
described by fractional differential equations.
In Section~\ref{sec8}, a kind of Sobolev type condition
and a nonlocal control condition for nonlinear fractional multiple control systems
is considered. The Sobolev condition is given in terms of two linear operators and
requires formulating two other characteristic solution operators and their properties,
such as boundedness and compactness. Further, we consider
an optimal control problem of multi-integral functionals, with integrands
that are not convex in the controls. An interrelation between
the solutions of the original problem and the relaxation one is given. 
Under certain assumptions, it is shown that the relaxed problem has a solution
with interesting convergence properties \cite{AMA.15,AMA.31,MyID:323}.


\section{Basic definitions and notations}

Here we review some essential facts from fractional calculus \cite{29,AMA.23},
basic definitions of a fuzzy number and fuzzy concepts \cite{24,26,27},
semigroup theory \cite{3-pa,AMA.44}, and multi-valued analysis \cite{AMA.3,AMA.21}.

\begin{definition}
\label{Definition 2.1}
The fractional integral of order $\alpha>0$ of a function $f\in L^{1}([a,b],\mathbb{R})$ is given by
$$
I^{\alpha}_{a}f(t):=\frac{1}{\Gamma(\alpha)}\int_{a}^{t}(t-s)^{\alpha-1}f(s)ds,
$$
where $\Gamma$ is the classical gamma function.
\end{definition}

If $a=0$, then we can write $I^{\alpha}f(t) := (g_{\alpha}*f)(t)$, where
$$
g_{\alpha}(t):=
\left\{
\begin{array}{ll}
\frac{1}{\Gamma(\alpha)}t^{\alpha-1},& \mbox{$t>0$},\\
0, & \mbox{$t\leq 0$}
\end{array}\right.
$$
and, as usual, $*$ denotes convolution. Moreover,
$\lim\limits_{\alpha \downarrow 0} g_{\alpha}(t)=\delta(t)$
with $\delta$ the delta Dirac function.

\begin{definition}
\label{Definition 2.2}
The Riemann--Liouville fractional derivative of order $\alpha>0$,
$n-1<\alpha<n$, $n\in \mathbb{N}$, is given by
$$
^{L}D^{\alpha}f(t) := \frac{1}{\Gamma(n-\alpha)}\frac{d^{n}}{dt^{n}}
\int_{0}^{t}\frac{f(s)}{(t-s)^{\alpha+1-n}}ds, \quad t>0,
$$
where function $f$ has absolutely continuous derivatives up to order $(n-1)$.
\end{definition}

\begin{definition}
\label{Definition 2.3}
The Caputo fractional derivative of order $\alpha>0$,
$n-1<\alpha<n$, $n\in \mathbb{N}$, is given by
$$
^{C}D^{\alpha}f(t) := ~^{L}D^{\alpha}\left(f(t)
-\sum\limits_{k=0}^{n-1}\frac{t^{k}}{k!}f^{(k)}(0)\right),
\quad t>0,
$$
where function $f$ has absolutely continuous derivatives up to order $(n-1)$.
\end{definition}

If $f$ is an abstract function with values in $X$, then the integrals
that appear in Definitions~\ref{Definition 2.1} to \ref{Definition 2.3}
are taken in Bochner's sense.

\begin{remark}
\label{Remark 2.1}
Let $n-1<\alpha<n$, $n\in \mathbb{N}$.
The following properties hold:
\begin{itemize}
\item [(i)] If $f\in C^{n}([0, \infty[)$, then
$$
^{C}D^{\alpha}f(t)=\frac{1}{\Gamma(n-\alpha)}
\int_{0}^{t}\frac{f^{(n)}(s)}{(t-s)^{\alpha+1-n}}ds
=I^{n-\alpha}f^{(n)}(t),
\quad t>0;
$$

\item [(ii)] The Caputo derivative of a constant function is equal to zero;
	
\item [(iii)] The Riemann--Liouville derivative of a constant function is given by
$$
^{L}D^{\alpha}_{a^{+}}C=\frac{C}{\Gamma(1-\alpha)}(t-a)^{-\alpha},~ 0<\alpha<1.
$$
\end{itemize}
\end{remark}

We denote the set of all real numbers by $\mathbb{R}$  
and the set of all fuzzy numbers on $\mathbb{R}$ 
is indicated by $\mathbb{E}$. A fuzzy number is a
mapping $u:\mathbb{R}\rightarrow[0,1]$ with the following
properties:\\
(i) ${u}$ is upper semi-continuous,\\
(ii) ${u}$ is fuzzy convex, i.e., ${u}(\lambda x+(1-\lambda)y\geq min \{{u}(x),{u}(y)\}$ 
for all $x,y \in \mathbb{R},\lambda\in [0,1],$\\
(iii) ${u}$ is normal, i.e., $\exists x_{0}\in \mathbb{R}$ for which ${u}(x_0)=1$,\\
(iv) $supp{u}= \{x \in \mathbb{R}|{u}(x)>0\}$ is the support of the $u$, and its closure
cl({\it supp ${u}$}) is compact.

\begin{definition}[See \cite{24}]
\label{d-2}
We define a metric $D$ on $\mathbb{E}$ 
($D:\mathbb{E}\times \mathbb{E}\longrightarrow \mathbb{R}_+\bigcup \{0\}$) 
by a distance, namely the Hausdorff distance as follows:
\begin{equation}
D({u},{v}) = \mathop {\sup }\limits_{r \in [0,1]} 
\;\max \{ |{{u}}_-(r) - {{v}}_-(r)|,| {{u}}_+(r) -  {{v}}_+(r)|\}
\end{equation}
It is shown that $(\mathbb{E},D )$ is a complete metric space.
\end{definition}

The concept of Hukuhara-difference, which is recalled in the next definition,  
was initially generalized by Markov \cite{41} to introduce the notion 
of generalized Hukuhara-differentiability for the interval-valued functions. 
Afterwards, Kaleva \cite{24} employed this notion to define the 
fuzzy Hukuhara-differentiability for the fuzzy-valued functions.  

\begin{definition}[See \cite{24}]
\label{d1}
Let $x,y\in \mathbb{E}$. If there exists $z\in \mathbb{E}$ such that 
$x=y\oplus z$, then $z$ is called the Hukuhara-difference of $x$ and $y$, 
and it is denoted by $x\ominus y$.
\end{definition}

\begin{definition}[See \cite{29,33}]
\label{d5}
We denote the Caputo fractional derivatives by the capital letter 
with upper-left index ${{}^cD}$, and the Caputo fractional derivatives 
of order $v$ is defined as
\begin{align*}
 {{}^cD^v}f(x) &= {I^{m - v}}{D^m}f(x) \\
 &= \frac{1}{{\Gamma (m - v)}} \int_0^x {{(x - t)}^{m - v - 1}}f^m(t)\,dt
\end{align*}
where $m - 1 < v \le m$; $x > 0$; and $D^m$ is the classical 
differential operator of order $m$.
\end{definition}

Let $a>0$ and $J=(0,a]$. We denote $C(J,\mathbb{E})$ as the space 
of all continuous fuzzy functions defined on $J$. Also let $f\in C(J,\mathbb{E})$. 
Then we say that $f \in {L^1}(J,\mathbb{E})$ if and only if 
$D(\int_0^a {f(s)~ds,\hat 0} ) < \infty $ \cite{40}. In the rest of the paper,  
the above notations will be used frequently. The fuzzy Caputo derivatives 
of order $0 < v \le 1$ for a fuzzy-valued function $f$ are given as follows.

\begin{definition}[See \cite{39}]
\label{d8}
Let $f\in C(J,\mathbb{E})\cap {L^1}(J,\mathbb{E})$ be a fuzzy set-value function. 
Then $f$ is said to be Caputo's fuzzy differentiable at $x$ when
\begin{equation}
({}^cD_{{0^ + }}^v f)(x) = \frac{1}
{{\Gamma (1 - v)}}\int_0^x {\frac{{f'(t)}}{{{{(x - t)}^v}}}} dt,
\label{eq5}
\end{equation}
where $0<v\le1$.
\end{definition}

We now proceed with some basic definitions and results from multivalued analysis.
For more details on multivalued analysis
we refer to the books \cite{AMA.3,AMA.21}. We use the following symbols:
$P_{f}(T)$ is the set of all nonempty closed
subsets of $T$; $P_{bf}(T)$ is the set of all nonempty, closed and bounded
subsets of $T$. On $P_{bf}(T)$, we have a metric, known as the Hausdorff metric,
defined by
\[
d_{H}(A,B):=\max\left\lbrace\sup_{a\in A}d(a,B),\,\sup_{b\in B}d(b,A)\right\rbrace,
\]
where $d(x,C)$ is the distance from a point $x$ to a set $C$.
We say that a multivalued map is $H$-continuous if it is continuous in
the Hausdorff metric $d_{H}(\cdot,\cdot)$.
Let $F: I \rightrightarrows 2^{T}\backslash\{\emptyset\}$ be a multifunction.
For $1\leq p\leq +\infty$, we define
$S^{p}_{F}:=\lbrace f\in L^{p}(I, T): f(t)\in F(t)$ a.e. on $I\rbrace$.
We say that a multivalued map $F:I \rightrightarrows P_f(T)$ is measurable if
$F^{-1}(E)=\{t\in I:F(t)\cap E\neq\emptyset\}\in \Sigma$ for every closed
set $E\subseteq T$. If $F:I\times T\to P_f(T)$, then the measurability
of $F$ means that $F^{-1}(E)\in\Sigma\otimes\mathcal{B}_{T}$,
where $\Sigma\otimes\mathcal{B}_{T}$ is the $\sigma$-algebra of subsets in
$I\times T$ generated by the sets $A\times B$, $A\in\Sigma$,
$B\in\mathcal{B}_{T}$, and $\mathcal{B}_{T}$ is the $\sigma$-algebra of
the Borel sets in $T$.

Suppose that $V_{1}$ and $V_{2}$ are two Hausdorff topological spaces and
$F: V_{1}\to 2^{V_{2}} \backslash\{\emptyset\}$. We say that $F$
is lower semicontinuous in the sense of Vietoris (l.s.c., for short)
at a point $x_0\in V_{1}$, if for any open set $W\subseteq V_{2}$,
$F(x_0)\cap W\neq\emptyset$, there is a neighborhood
$O(x_0)$ of $x_0$ such that $F(x)\cap W\neq\emptyset$
for all $x\in O(x_0)$. Similarly, $F$ is said to be upper semicontinuous in the sense
of Vietoris (u.s.c., for short) at a point $x_0\in V_{1}$, if for any open set
$W\subseteq V_{2}$, $F(x_0)\subseteq W$, there is a neighborhood $O(x_0)$ of $x_0$
such that $F(x)\subseteq W$ for all $x\in O(x_0)$. For more properties of
l.s.c and u.s.c, we refer to the book \cite{AMA.21}.
Besides the standard norm on $L^q(I, T)$ (here, $T$ is a separable reflexive
Banach space), $1<q<\infty$, we also consider the so called weak norm:
\begin{equation}
\label{eq:2.2}
\|u_{i}(\cdot)\|_{\omega}:=\sup_{0\leq t_1\leq t_2\leq a}
\left\Vert\int_{t_1}^{t_2}u_{i}(s)ds\right\Vert_T,~  u_{i}\in L^q(I, T),
\ i=1,\ldots,r.
\end{equation}
The space $L^q(I, T)$ furnished with this norm will be denoted by
$L_{\omega}^q(I, T)$. 


\section{Fuzzy fractional differential equations}
\label{int} 

Let
\begin{equation}
\label{eq21} 
{}_0^cD_t^\alpha X(t) = AX(t) + f(t)
\end{equation}
with initial condition $X(0)=X_0$ and the matrix $A$
\[
A = \left[ 
{\begin{array}{*{20}{c}}
{ - \frac{1}{{{s_2}}}}&{{s_0}}&0\\
{ - {s_0}}&{ - \frac{1}{{{s_2}}}}&0\\
0&0&{ - \frac{1}{{{s_1}}}}
\end{array}} 
\right] \in {\mathbb{R}^{3 \times 3}}.
\]
Also, let $f(t) = {\left(0,0,\frac{{{X_0}}}{{{s_1}}}\right)^T},\;X(t) 
={({X_x}(t),{X_y}(t),{X_z}(t))^T}$ and 
${X_0} = {({X_x}(0),{X_y}(0),{X_z}(0))^T}$ be fuzzy vectors.
Note that the coefficients of $A$  are expressed as follows:
\[
{s_0} = \frac{{{w_0}}}{{\tau _2^{\alpha  - 1}}},\;
\frac{1}{{{s_1}}} = \frac{{\tau _1^{\alpha  - 1}}}
{{{T_1}}},\frac{1}{{{s_2}}} 
= \frac{{\tau _2^{\alpha-1}}}{{{T_2}}},
\quad \alpha  \in (0,1]
\]
where $\tau_1$ and $\tau_2$ are fractional time constants.
The predictor-corrector method is  investigated in \cite{121}.
Particularly, the fractional Adams--Bashforth as a predictor and the
fractional Adams--Moulton as a corrector are exploited. Moreover,
a new variant of the fuzzy fractional Adams--Bashforth--Moulton
(FFABM) method is  introduced  \cite{121}. Finally, \cite{121} demonstrates 
the capability of the developed numerical methods 
for fuzzy fractional-order problemscin terms of accuracy and stability analysis.

\begin{figure}
\centering 
\includegraphics[width=3in]{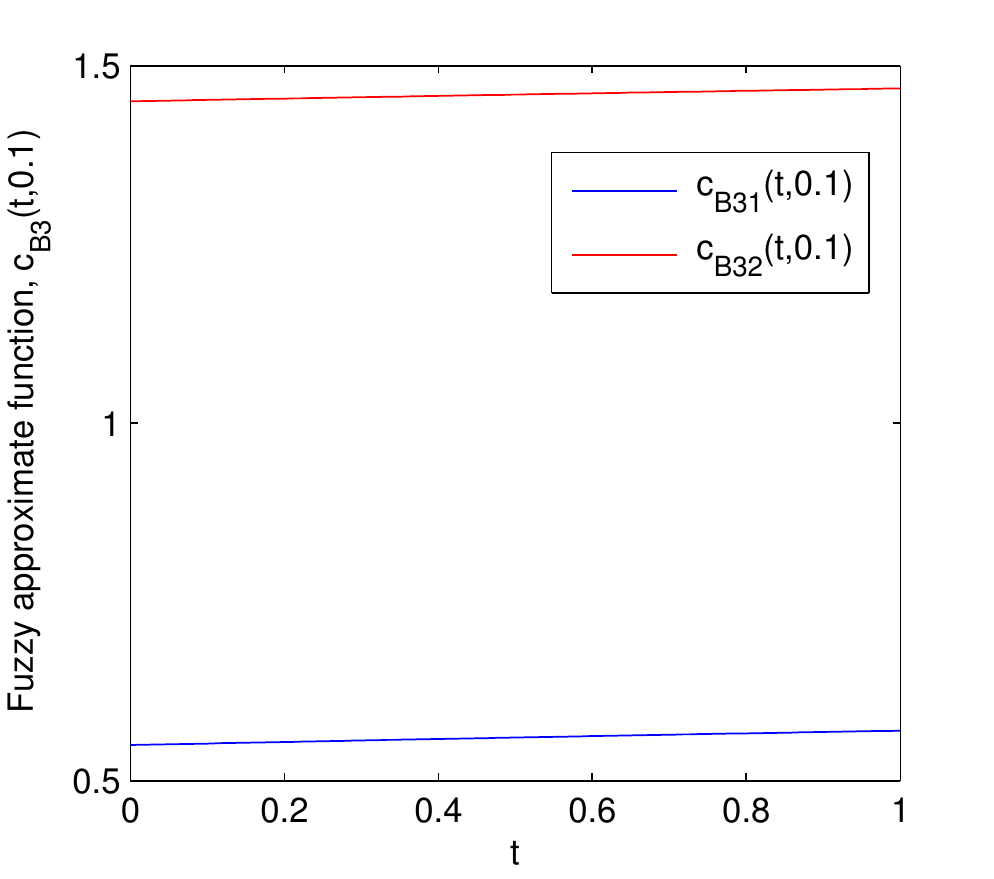}
\caption{\small{Case I: Fuzzy approximate function, 
${{{c_B}_3}}(t,0.1)$, over  $t\in[0,1]$ with
$N=3$, $v=0.85$ and $\alpha=0$, \cite{14}.}} 
\label{fig1}
\end{figure}

\begin{figure}
\centering
\includegraphics[width=3in]{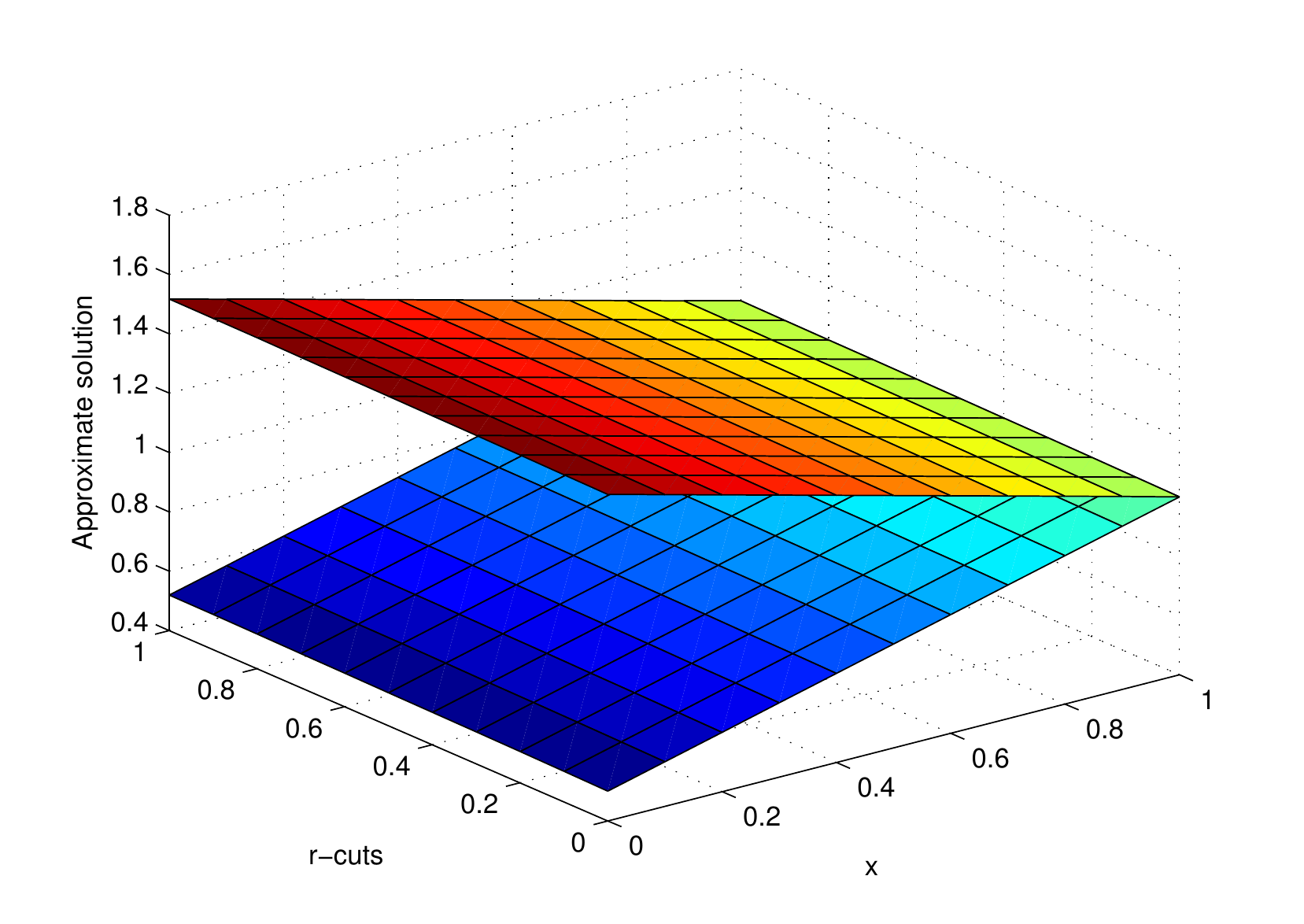}
\caption{\small{Case I:  Fuzzy approximate solution,
with $N=8$, $\alpha=0$ and $v=0.85$, \cite{14}.}} 
\label{fig2}
\end{figure}

\begin{figure}
\centering 
\includegraphics[width=3in]{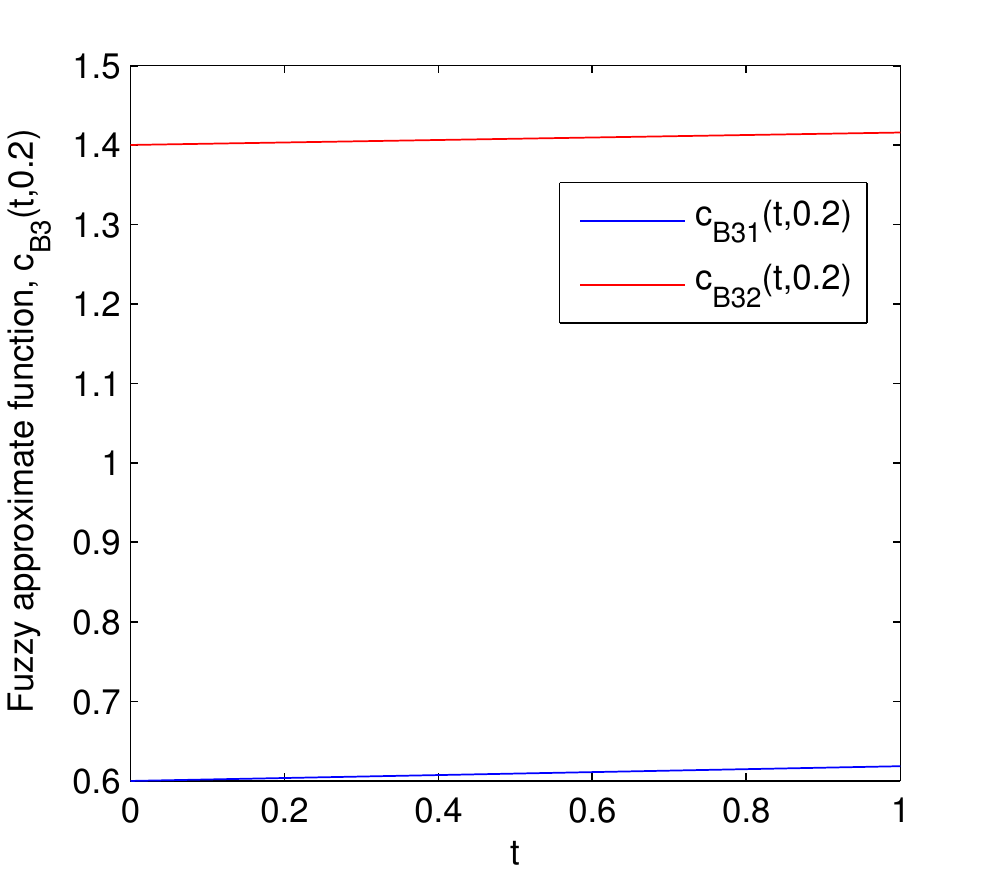}
\caption{\small{Case II: Fuzzy approximate function, 
${{{c_B}_3}}(t,0.2)$, over  $t\in[0,1]$ with
$N=3$, $v=0.95$ and $\alpha=0$, \cite{14}.}}
\label{fig3}
\end{figure}

Ahmadian et al. \cite{14} have dealt with the application of FFDEs 
to model and analyze a kinetic model of diluted acid hydrolysis 
under uncertainty as follows. When water is added to the Hemicellulose xylane,  
Xylose is formed through the hydrolysis reaction. Furfural is the main degradation
product obtained through the degradation of a molecule of Xylose
by the releasing of three water molecules. Scheme~\ref{eq2} demonstrates
the depolymerization of Xylan. The simplest kinetic model for the
hemicellulose hydrolysis was firstly proposed by Saeman
\cite{saeman}. He discovered that a straightforward two-step
reaction model sufficiently  explained the generation of sugars
during wood hydrolysis. Saeman's model assumed pseudo homogeneous
irreversible first-order consecutive reactions:
\begin{equation}
\label{eq2}
\text{Hemicellulose  Xylan} ~\xrightarrow{{{\text{Hydrolysis}}}}
~\text{Xylose} ~\xrightarrow{{{{-3H_2O}}}}~ {\text Furfural}.
\end{equation}
In real problems, we firstly choose the initial conditions as starting points.
Indeed, initial conditions for such models are determined by analyzer systems,  
which are not adequate for high accuracy results. So, instead of using 
deterministic values, it is better to  employ uncertain conditions. 
In order to consider the original problem in a new sense, the authors used
the fuzzy initial value instead of the crisp initial value. 
In this direction, they reconstructed the original problem based on
fuzzy fractional calculus:
\begin{equation}
\label{eq10}
{{}^cD^v}C_B(t) \oplus {k_2}{C_B(t)} 
= {k_1}{C_{{A_0}}}\odot\exp ( - {k_1}t), 
\end{equation}
where $C_B:{L^{\mathbb{R}_\mathcal{F}}}[0,1] \cap
{C^{\mathbb{R}_\mathcal{F}}}[0,1]$ is a continuous fuzzy-valued
function and ${}^cD_{{0^ + }}^v$ indicates the fuzzy Caputo
fractional derivative of order $v \in [0,1]$. Also, the Xylose
concentration is defined  at time $t=0$ as a fuzzy
number, i.e., $C_B(0;r) = [{C_{{B_{01}}}(r)},{C_{{B_{02}}}(r)}]$, 
$0 < r \le 1$. Additionally, the notations  $\oplus,\odot$ mean 
addition and multiplication, respectively in $\mathbb{R}_\mathcal{F}$.
This approach can also be a feasible alternative to ordinary differential 
equation models under uncertainty. As a case study, they have formulated 
a mathematical model to analyze the Xylose concentration in the acid 
hydrolysis as a promising source of Xylose. The proposed model has been 
described by a FFDE of order $0 < v \leq 1$ with the aim of incorporating 
uncertainty into the initial values and function of the model. 
This model provides a more realistic view taking into account 
the variations of Xylose concentration. Moreover, it constitutes 
a new approach for modeling a chemical reaction. This approach 
is used to obtain reliable data of acid concentration immediately 
during the process (see, Figs.~\ref{fig1}, \ref{fig2} and \ref{fig3}). 
Therefore, Ahmadian et al. \cite{15}  proposed a new iterative method  
for solving FFDEs of Caputo type. The basic idea is to convert FFDEs 
to a type of fuzzy Volterra integral equation. Then the obtained Volterra 
integral equation is exploited with some suitable quadrature rule 
to get a fractional predictor-corrector method.

Salahshour et al. \cite{17} investigated the solution of a class of
fuzzy sequential fractional differential equation of order $\beta \in
(0; 1)$. For this purpose, some basic results were developed
in terms of the fuzzy concept. The main scope of this study
is to present the solution of fuzzy problems based on the
contraction principle in the new fuzzy complete metric space. 
Afterwards, Ahmadian et al.\cite{18} conducted a research 
that has been devoted to solve linear FFDEs of
Caputo sense. The basic idea is to develop a fractional linear
multistep method for solving linear FFDEs under fuzzy fractional
generalized differentiability. The authors proposed fractional linear
multistep methods and investigated the consistence,
convergence and stability properties of the method. 
A new fractional derivative, with some simplifications 
in the formula and computations, was proposed under interval uncertainty 
in \cite{19}. The authors have developed the highlights of the conformable 
fractional derivative, which are more influential for the solution 
of fractional interval differential equations (FIDEs) under generalized 
Hukuhara differentiability. The reader interested on conformable 
fractional differentiation is refereed to \cite{MyID:324,MyID:351} 
and references therein. In practice, a simple fractional derivative 
satisfying the main rules for uncertain fractional differential equations, 
such as the product rule and the chain rule, simplifies considerably 
the cumbersome mathematical expressions of the mathematical modeling 
in engineering sciences. At the same time, several authors developed 
theoretical aspects of this concept for the solution of FDEs
under uncertainty in their researches \cite{20,21,22,23}, which were
a continuation of the concepts proposed in the monographs and papers 
of fuzzy setting theory \cite{24,25,26,27,28,281,282}  
and fractional calculus \cite{29, 30,31,32,33,34,35}.

Very recently, Salahshour et al. \cite{38} made a meaningful 
contribution to study a theory of FDEs under interval uncertainty. 
The new interval fractional derivative has several impressive properties 
that each of them lonely can improve significantly the modeling of 
real life systems using FDEs such as non-singularity kernel. 
It was devoted to research on this new and applicable notion 
in terms of interval uncertainty. The authors proposed 
adequate conditions for the uniqueness of solution 
of fractional interval differential equations (FIDEs) 
under Caputo--Fabrizio fractional derivative. 
Some interesting properties of this derivative were studied 
to make the easiest way for dealing with mathematical models 
based on Caputo--Fabrizio fractional interval equations. 
Ahmadian et al. \cite{39}  presented numerical simulations 
and introduced fuzzy mathematical models that can be represented in
terms of FDEs under certainty. Consider the following linear FFDE:
\begin{equation}
\left\{ {\begin{array}{*{20}{c}}
{({}^cD_{{0^ + }}^v  y)(x) + y(x) =
 f(x),\quad 0 < v  \le 1,}\\
{y(0) = {{ y}_0} \in \mathbb{E},}
\end{array}} \right.
\label{eq20}
\end{equation}
in which ${y}\in C(J,\mathbb{E})\cap {L^1}(J,\mathbb{E})$ 
is a continuous fuzzy-valued function, ${}^cD_{{0^ + }}^v$ 
indicates the fuzzy Caputo's fractional derivative 
of order $v$ and $f(x):[0,1] \mapsto \mathbb{E}$.
Compared with the extensive amount of work put into developing 
FDE schemes in the literature, we found out that only a little 
effort has been put into developing numerical methods for FFDE. 
Even so, most of the solutions are based on a rigorous framework, 
that is, they are often tailored to deal with specific applications 
and are generally intended for small-scale fuzzy fractional systems. 
In this research, the authors deployed a spectral tau method based 
on Chebyshev functions to reduce the FFDE to a fuzzy algebraic linear equation
system to address fuzzy fractional systems. The main advantage of this  
technique, using shifted Chebyshev polynomials in the interval $[0,1]$, 
is that only a small number of the shifted Chebyshev polynomials is 
required as well as the good accuracy that will be acquired 
in one time program running. Thus, it greatly simplifies 
the problem and reduces the computational costs. The solution is expressed 
as a truncated Chebyshev series and so it can be easily evaluated 
for arbitrary values of time using any computer program without 
any computational effort. The algorithm of the technique is as follows:
we generate $N$ fuzzy linear equations by applying
\begin{equation}
{\left\langle {R_N(x,r),T_i^*(x)} \right\rangle _E} 
= \tilde 0,\quad i = 0,1,\ldots,N - 1,\;r \in [0,1],
\label{eq221}
\end{equation}
where ${\left\langle {R_N(x,r),T_i^*(x)} \right\rangle _E} 
= {[(FR)\int_0^1 {{R_N}(x,r) \odot T_i^*(x)\odot w(x)dx} ]}$,  
$T_i^*(x)$ is a shifted Chebyshev polynomial and $R_N$ 
is a \emph{fuzzy-like residual operator} for \eqref{eq20}, 
which is defined in the matrix operator form
\[
{R_N}(x,r) = [{\underline{R}_N}(x,r),{\overline{R}_N}(x,r)],
\]
where
\begin{equation}
\left\{ {\begin{array}{*{20}{c}}
{{\underline{R}_N}(x,r) = {\underline{C}^T}(r)({D^{(v)}}\Phi (x) 
+ \Phi (x)) - {\underline{F}^T}(r)\Phi (x),}\\
{{\overline{R}_N}(x,r) = {\overline{C}^T}(r)({D^{(v)}}\Phi (x) 
+ \Phi (x)) - {\overline{F}^T}(r)\Phi (x).}
\end{array}} \right.
\label{eq88}
\end{equation} 
Using the definition of fuzzy-like inner product, we have:
\begin{equation}
{\left\langle {{D^{(v)}}{y_N}(x,r),T_k^*(x)} \right\rangle _E} 
+ {\left\langle {{y_N}(x,r),T_k^*(x)} \right\rangle _E} 
= {\left\langle {f(x,r),T_k^*(x)} \right\rangle _E}\quad k = 0,1,\ldots,N - 1,
\label{eq70}
\end{equation}
and $r\in[0,1]$. Then, in order to acquire the approximation 
$y_N(x,r)$ using the shifted Chebyshev tau approximation,  
we should find the unknown vector $C^T=[\underline{C}^T(r),{\overline{C}^T}(r)]$. 
Therefore, \eqref{eq70} can be stated as follows:
\begin{equation}
\left\{ {\begin{array}{*{20}{c}}
\begin{array}{l}
\sum\limits_{j = 0}^N {{c_j}(r)\left[ {{{\left\langle {{D^{(v)}}T_j^*(x),T_k^*(x)} \right\rangle}_w} 
+ {{\left\langle {T_j^*(x),T_k^*(x)} \right\rangle}_w}} \right]} \\
= {\left\langle {f(x,r),T_k^*(x)} \right\rangle_w},
\quad k = 0,1,\ldots,N - 1,\;j = 0,1,\ldots,N,\;r \in [0,1],
\end{array}\\
\hspace{-5.7cm}{\sum\limits_{j = 0}^N {{c_j}(r)T_j^*(0) = {y_0}(r).} }
\end{array}} \right.
\label{eq71}
\end{equation}
Then, using the the matrix form and their defined elements, 
\eqref{eq71} can be written in the following matrix form:
\begin{equation}
(\mathfrak{A} + \mu \mathfrak{B})\boldsymbol{C} = \boldsymbol{f}.
\label{eq72}
\end{equation}
Finally, system \eqref{eq72} can be solved based on the following lower-upper 
representation by any direct or numerical method:
\[
\left\{ {\begin{array}{*{20}{c}}
{\underline{(\mathfrak{A} +  \mathfrak{B})\boldsymbol{C}} = \underline{\boldsymbol{f}},}\\
{\overline{(\mathfrak{A} +  \mathfrak{B})\boldsymbol{C}} = \overline{\boldsymbol{f}}.}
\end{array}} \right.
\]
The analysis of the behaviors of physical phenomena is important  
in order to discover significant features of the character and 
structure of the  mathematical models. In a very recent time,  
Ahmadian et al. \cite{42} defined a new fuzzy
approximate solution and fuzzy approximate functions
formed on the generalized fractional Legendre polynomials (GFLPs) 
introduced in \cite{chen}, and then fuzzy Caputo fractional-order 
derivatives of GFLPs in terms of GFLPs themselves are stated and proved. 
They derived an effective spectral tau method under uncertainty
by applying these functions to solve two important fractional dynamical models  
via the fuzzy Caputo-type fractional derivative. They proposed  
a new model based on fractional calculus to deal with the Kelvin-Voigt (KV) 
equation and non-Newtonian fluid behavior model with fuzzy parameters
(for Caputo fractional Voigt models see \cite{MyID:350}).
Numerical simulations are carried out and the analysis of the results
highlights the significant features of the new technique 
in comparison with the previous findings.
The homogeneous strain relaxation equation with memory can be
defined as a differential equation of fractional order under
uncertainty as follows:
\begin{equation}
\label{eq11}
({}^cD_{{0^ + }}^v \omega)(t) \oplus B\odot\omega (t) 
= \frac{1}{\beta }\odot\rho (t),
\end{equation}
where $\omega:{L^{\mathcal{K}}}[0,h] \cap
{C^{\mathcal{K}}}[0,h]$ is a continuous fuzzy-valued function,
${}^cD_{{0^ + }}^v$ represents the Caputo-type fractional
derivative under fuzzy notion with  order $v \in [0,1]$ and $h\in
(0,1]$. Since the \eqref{eq11} is a general model of viscosity
behavior for non-Newtonian fluid under uncertainty, it is more
convenient to define $\beta\neq 0$  and $\rho(t)$ as a fuzzy number
and fuzzy set-valued function, respectively. Therefore, $B$ can be
a fuzzy number that alter the conditions of the model.
In addition, they developed the fuzzy fractional KV model under the
Caputo gH-differentiability that offers fuzzy models for
mathematical systems of natural phenomena as follows:
\begin{equation}
\label{eq14}
\rho(t) = {E}\odot \delta(t) \oplus {\eta}D_{{0^ + }}^v\delta(t)
\end{equation}
in which $\delta:{L^{\mathcal{K}}}[0,h] \cap
{C^{\mathcal{K}}}[0,h]$ presents a continuous fuzzy function,
${}^cD_{{0^ + }}^v$ specifies the Caputo-type derivative  
for $v \in [0,1]$ and $h\in (0,1]$. Also, $E $ and $\rho (t)$ 
can be defined as a fuzzy parameter and a fuzzy set-valued function, respectively, 
based on the conditions of the model. To test the proposed technique, practically, 
they solved the following two problems: 
\begin{equation}
\label{eq34}
\left\{ {\begin{array}{*{20}{c}}
{({}^cD_{{0^ + }}^v\omega )(t) \oplus  \omega (t) =  \rho (t),\quad t\in[0,1],}\\
{\omega (0,r) = [ - 1 + r,1 - r],\quad r\in[0,1],}
\end{array}} \right.
\end{equation}
where $\omega(t)$ is the fuzzy stress-strain function. In addition,
the generalized viscous coefficient and the modulus of
elasticity are assumed to be one ($E=\beta=1)$. Also, at first, they
considered that $\rho(t)=te^{-t}$. The next model is as follows:
\begin{equation}
\label{eq41}
\left\{ {\begin{array}{*{20}{c}}
\sqrt{t}+\frac{\sqrt{\pi}}{2} =  \delta(t) \oplus {}^cD_{{0^ + }}^v\delta(t),
\quad\; t\in[0,1],\\
{\delta (0,r) = [ -1+r,1-r],\quad r\in[0,1],}
\end{array}} \right.
\end{equation}
in which it is assumed that ${\eta}=1$, ${E}=1$ 
and $\rho(t)=\sqrt{t}+\frac{\sqrt{\pi}}{2}$. 
\begin{figure}
\centering        
\scalebox{.5}  
{\includegraphics{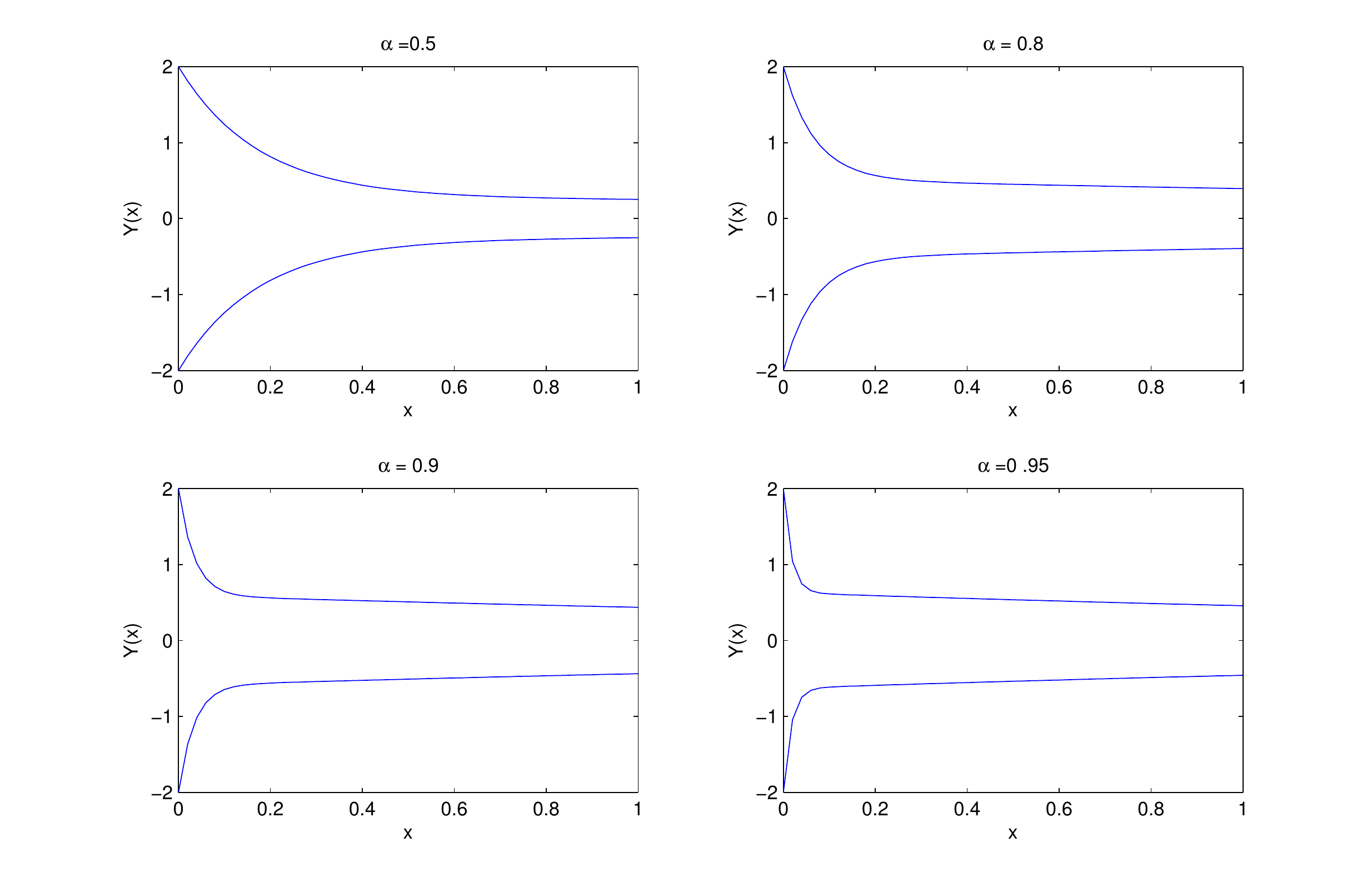}}
\caption{\small{Exact solution of the Basset problem based 
on the FIDE \cite{38}. }} 
\label{fig6}
\end{figure}
\begin{figure}
\centering
{\includegraphics[height=0.2\linewidth, width=0.5\linewidth]{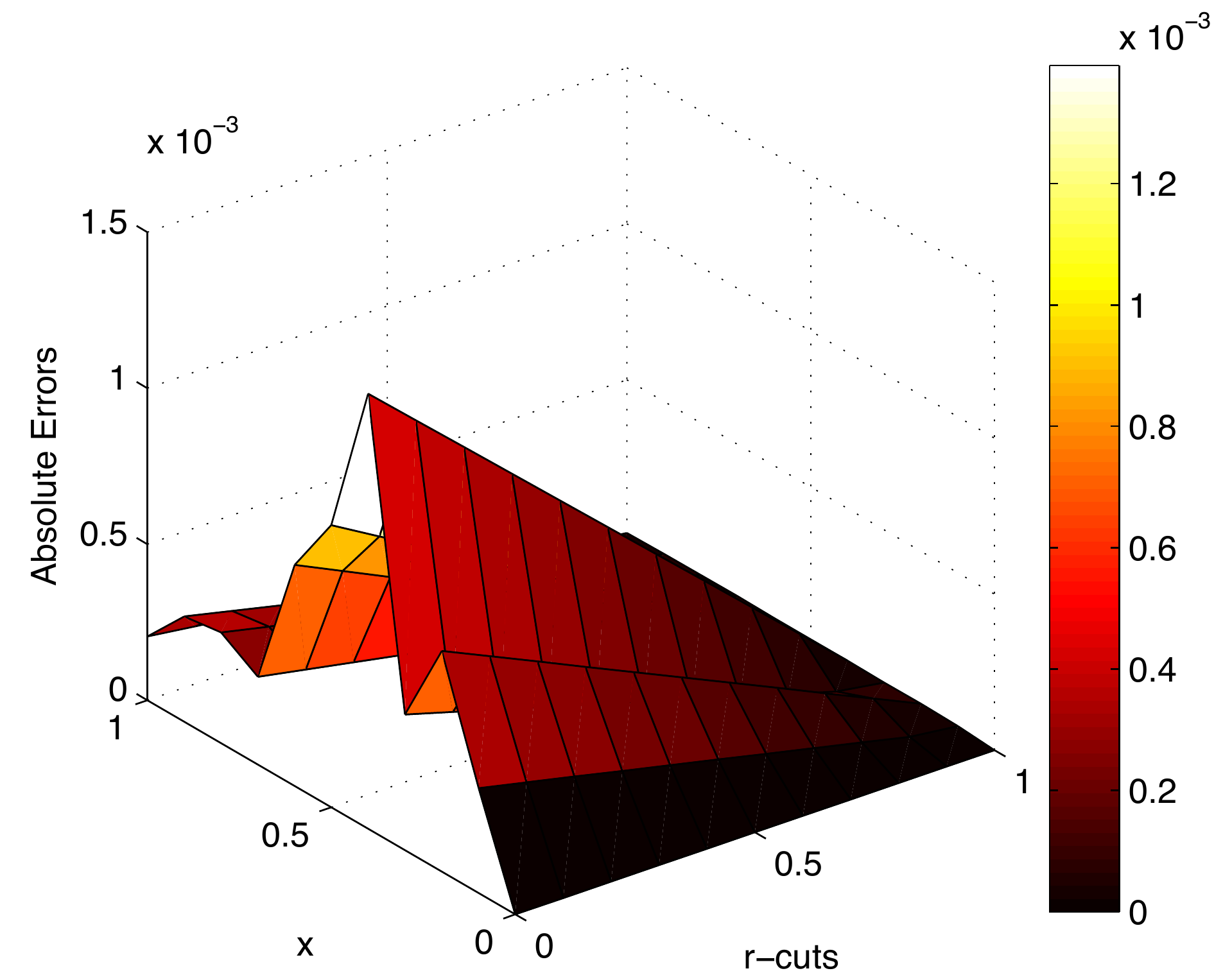}}
\caption{\small{Viscoelasticity: ${\underline{E}}_c(t;r)$  for $r\in[0,1]$ 
and $t\in[0,1]$ with $v=0.85$ and $N=10$ \cite{39}.}}
\label{fig4}
\end{figure}
\begin{figure}
\centering     
\subfigure[$N=10$]{\includegraphics[height=0.3\linewidth, width=0.4\linewidth]{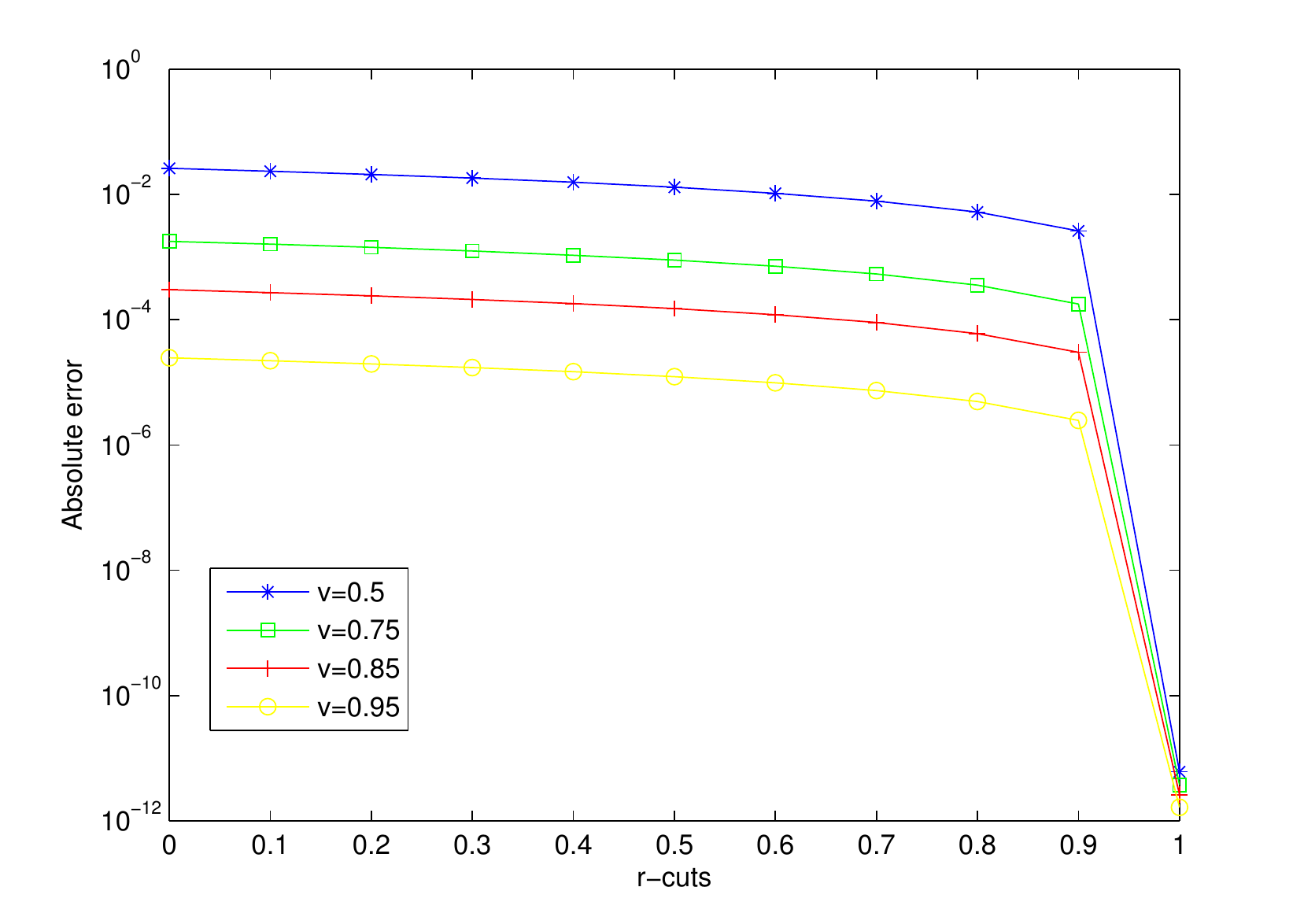}
\label{fig11}}
\subfigure[$v=0.95$]{\includegraphics[height=0.3\linewidth, width=0.4\linewidth]{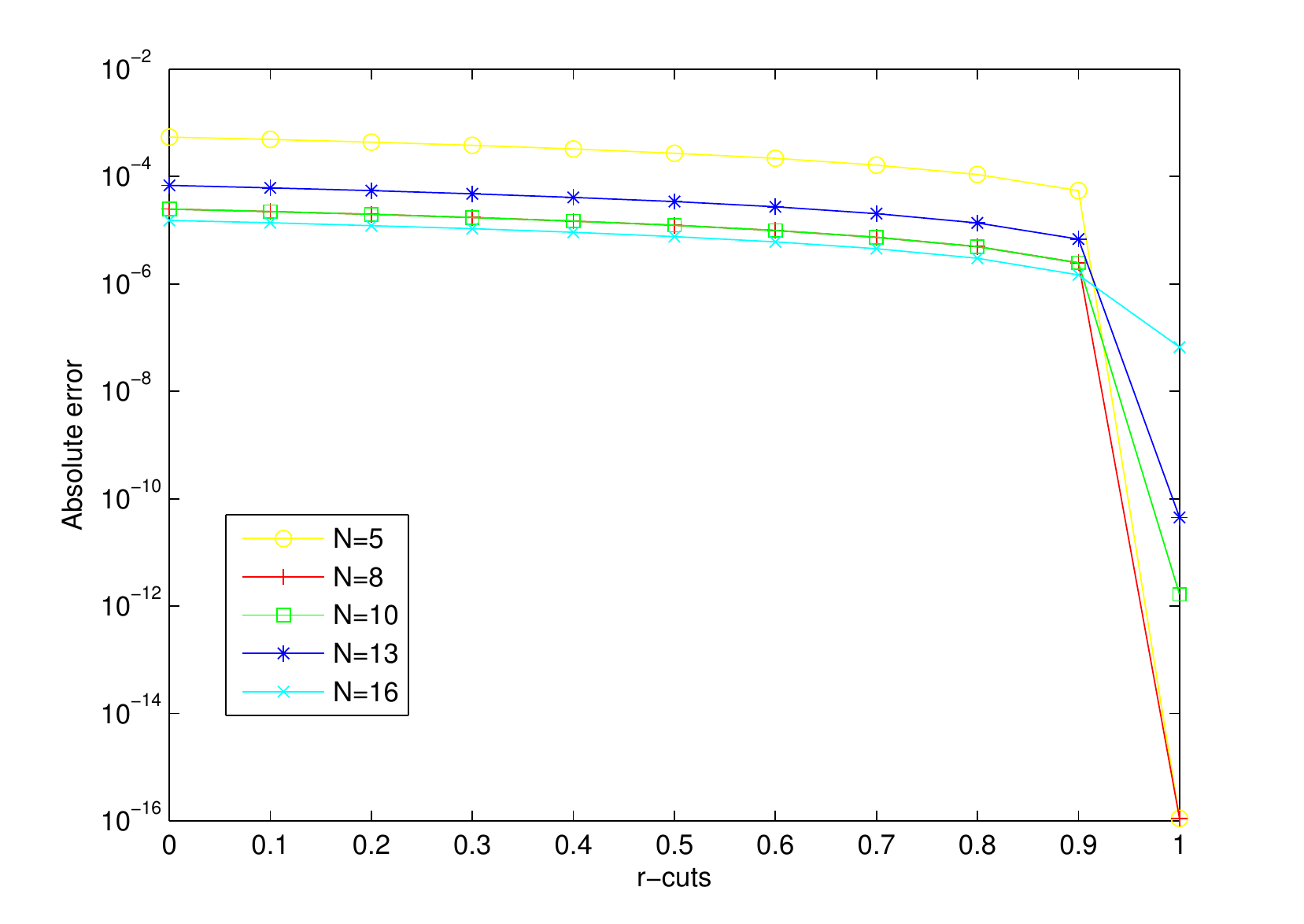}
\label{fig21}} 
\caption{\small{Viscoelasticity: Absolute errors of the proposed method, 
$\underline{E}_c(1;r)$ (a) for different values of Caputo derivative $v$; 
and (b) for different values $N$ \cite{39}.}}
\label{fig2a}
\end{figure}
\begin{figure*}
\centering 
\subfigure[]{\includegraphics[height=0.2\linewidth, width=0.4\linewidth]{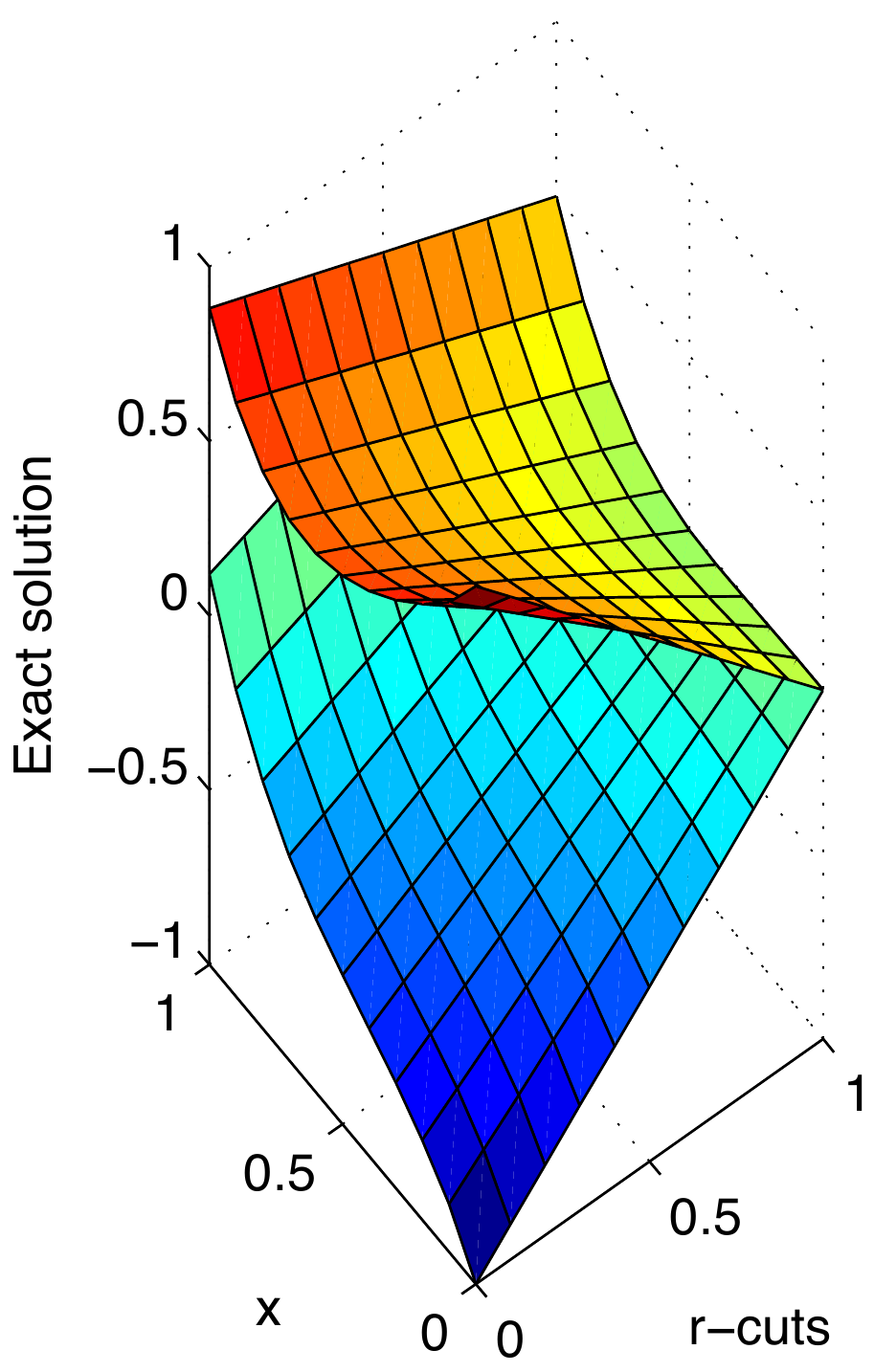}}
\subfigure[]{\includegraphics[height=0.2\linewidth, width=0.4\linewidth]{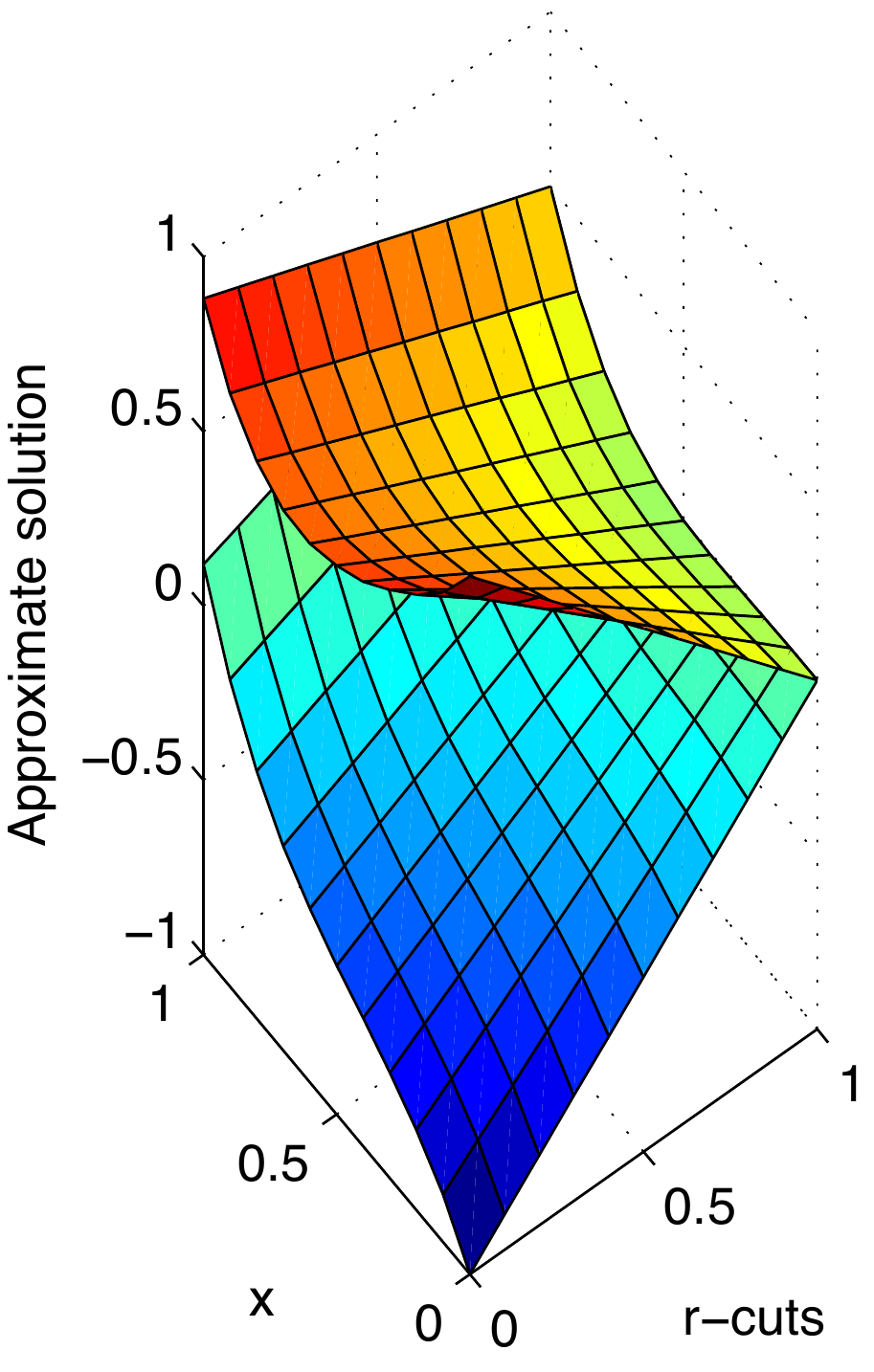}}
\caption{\small{Profiles of (a) the exact solution and 
(b) the fuzzy solution of the motion model of a rigid plate immersed 
in a Newtonian fluid with $v=0.85$ and $N=8$ \cite{39}.}}
\label{fig5}
\end{figure*}
\begin{figure}
\centering         
\scalebox{.25} 
{\includegraphics{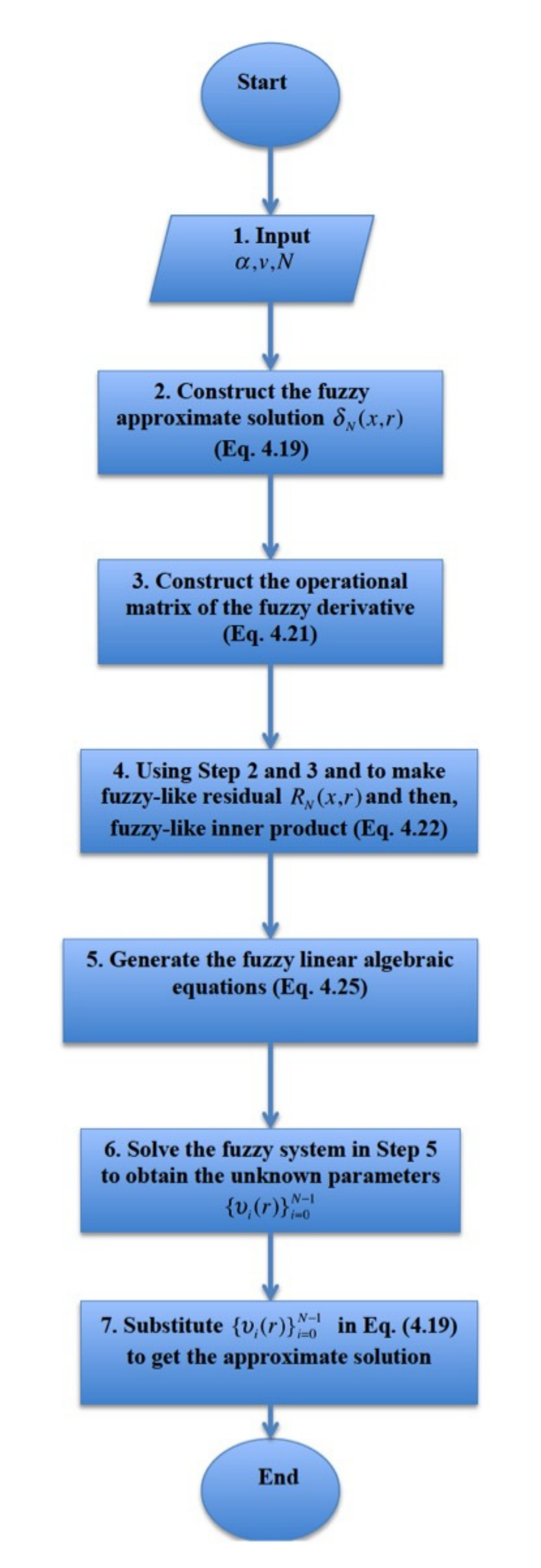}}
\caption{\small{Flowchart of the implementation of the proposed technique \cite{42}.}}
\label{fig31}
\end{figure}
\begin{figure}[!ht]
\centering  
\subfigure[]{\includegraphics[height=0.35\linewidth, width=0.45\linewidth]{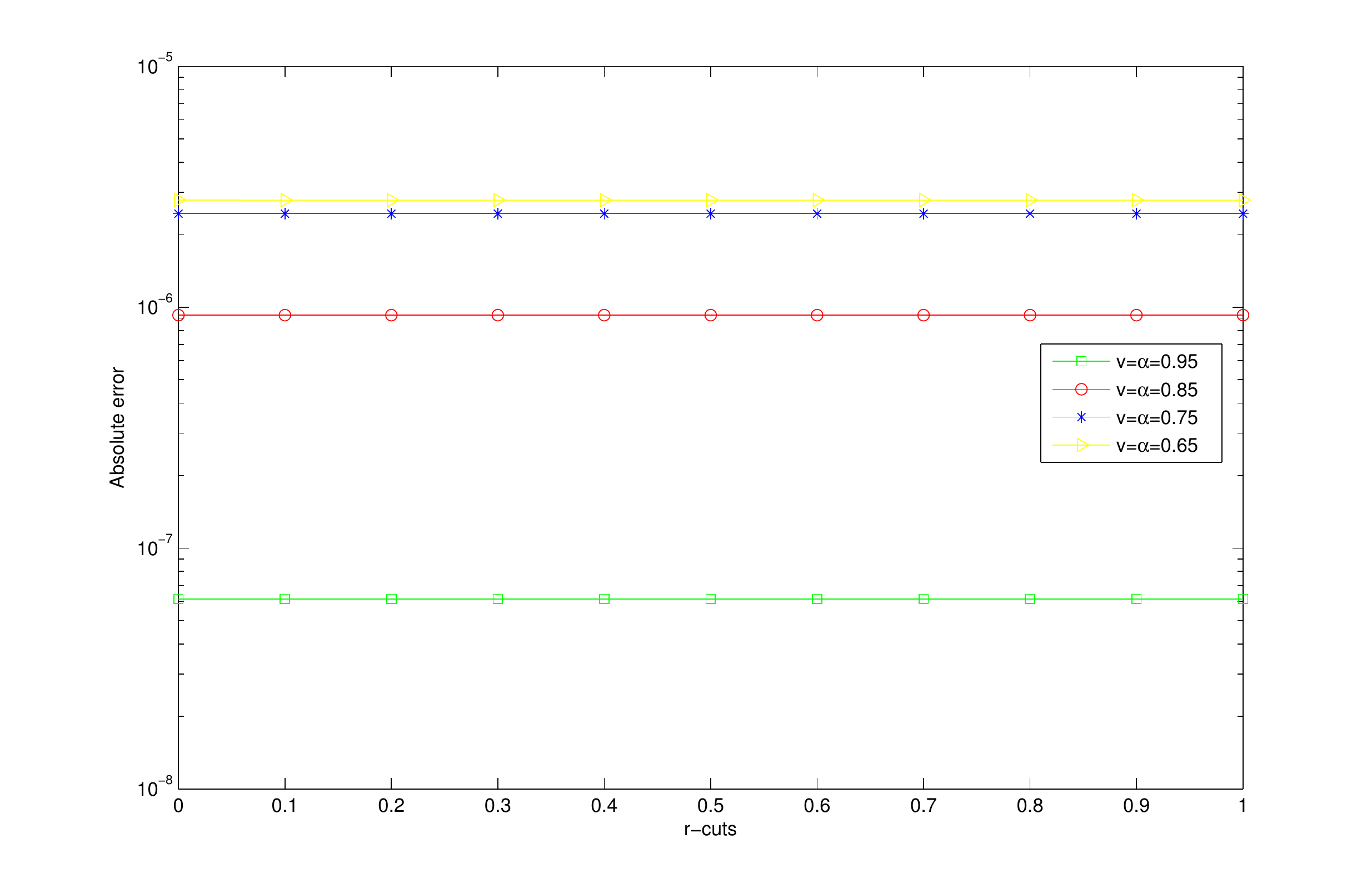}}
\subfigure[]{\includegraphics[height=0.35\linewidth, width=0.45\linewidth]{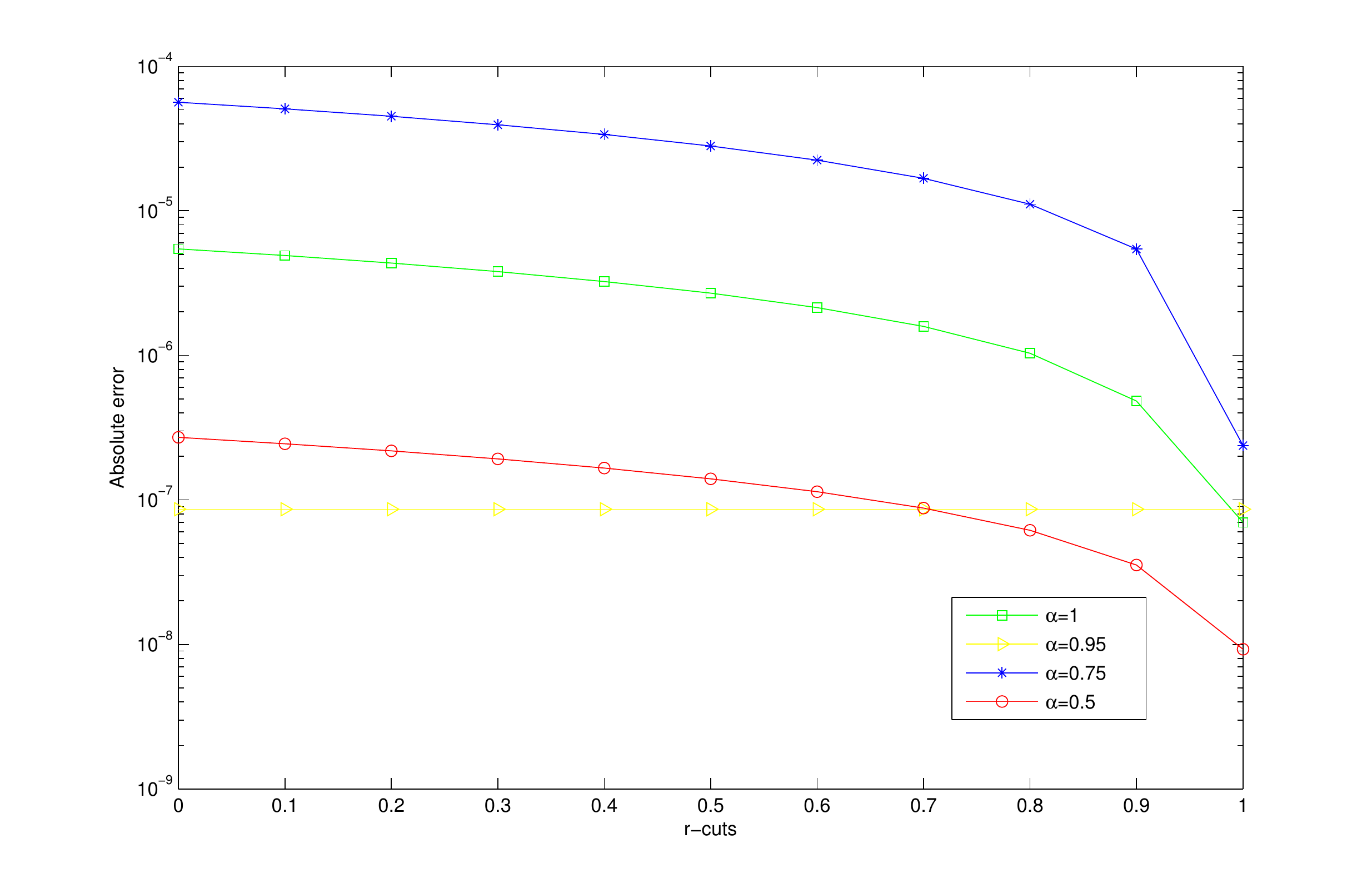}}
\caption{\small{Comparison of the absolute error functions,  
$\underline{\mathtt{E}_g}(1,r)$ : (a) $v$ and  $\alpha$ have the same values 
and $N=8$~ (b) $v=0.95$, $N=9$ for different fractional orders of GFLFs, 
at $T=1$ for the problem \eqref{eq34} \cite{42}.}}
\label{fig51}
\end{figure}
\begin{figure}[!ht]
\subfigure[  $v=\alpha=0.75$
\label{subfig7a}]{\includegraphics[width=0.45\textwidth]{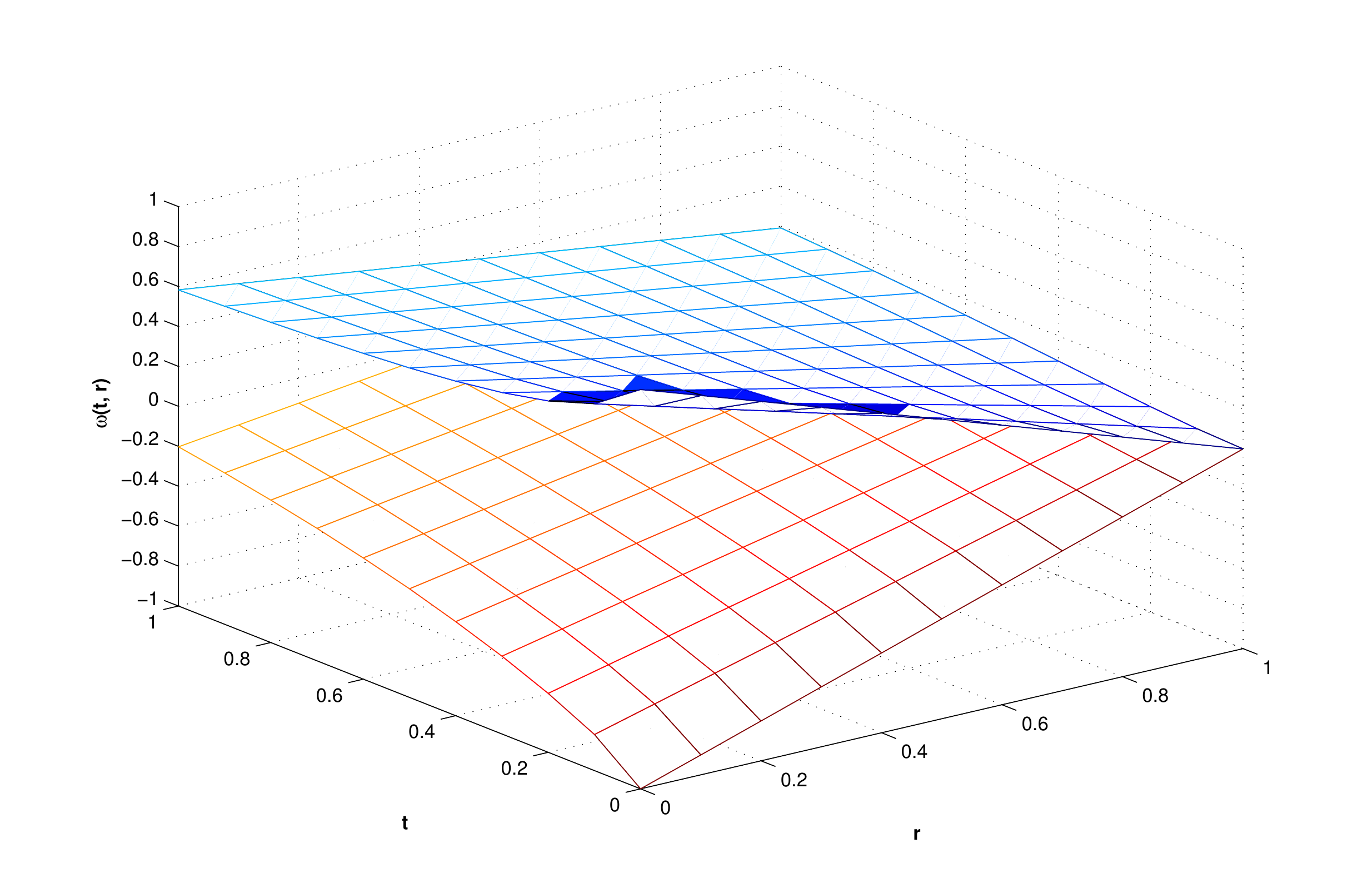}}
\hfill
\subfigure[$v=\alpha=0.95$ \label{subfig7b}]{%
\includegraphics[ width=0.45\linewidth]{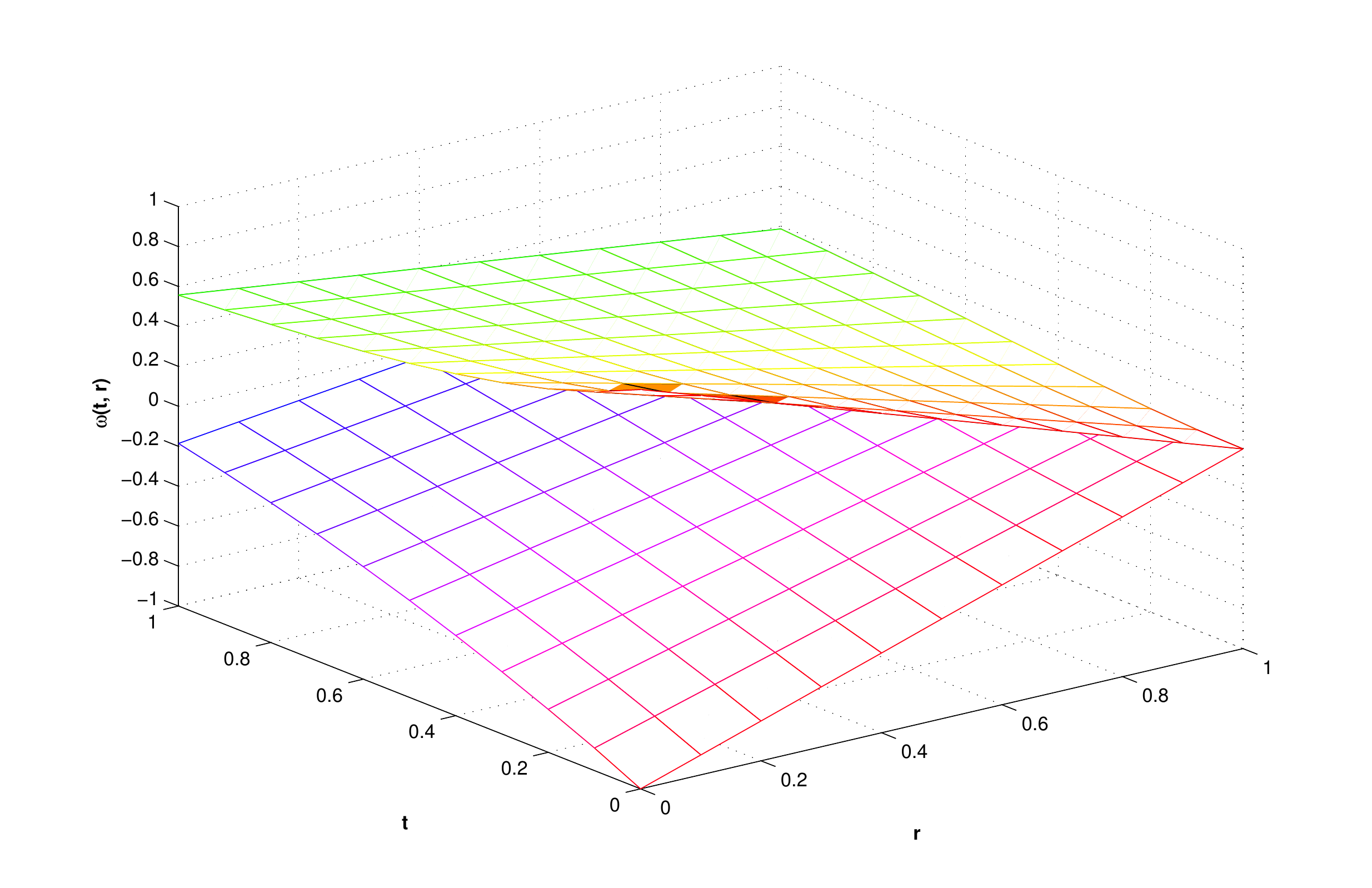}
}
\caption{\small{Fuzzy approximate solutions of the problem \eqref{eq34} 
for two different orders of the fractional derivative with 
$N=8$ over $t\in[0,1]$ and $r\in [0,1]$\cite{42}.}}
\label{fig7}
\end{figure}
\begin{figure}[!ht]
\subfigure[  $\underline{\mathtt{E}_g}(t,r)$
\label{subfig11a}]{%
\includegraphics[width=0.45\textwidth]{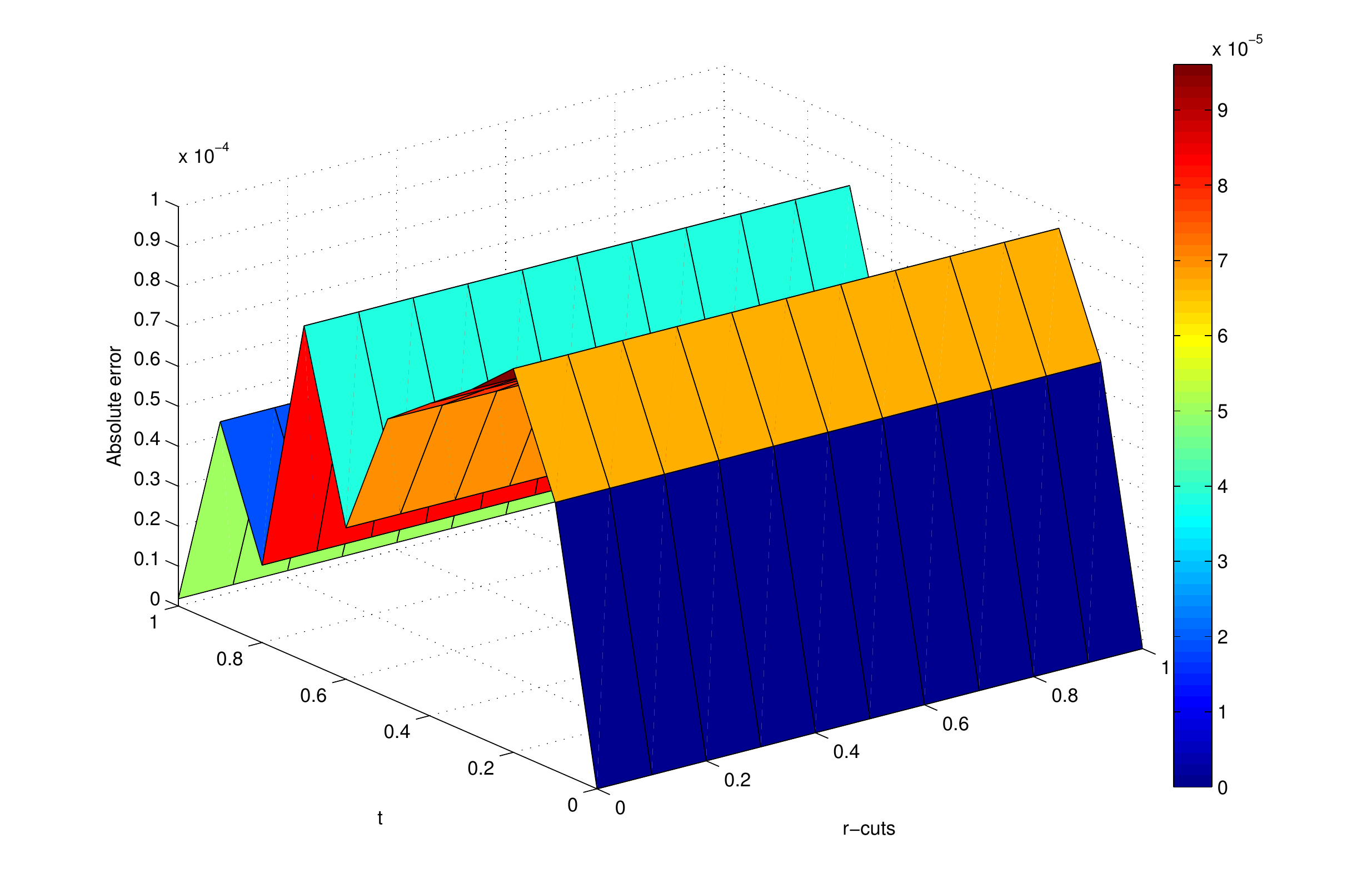}
}
\hfill
\subfigure[$\overline{\mathtt{E}_g}(t,r)$  \label{subfig11b}]{%
\includegraphics[ width=0.45\linewidth]{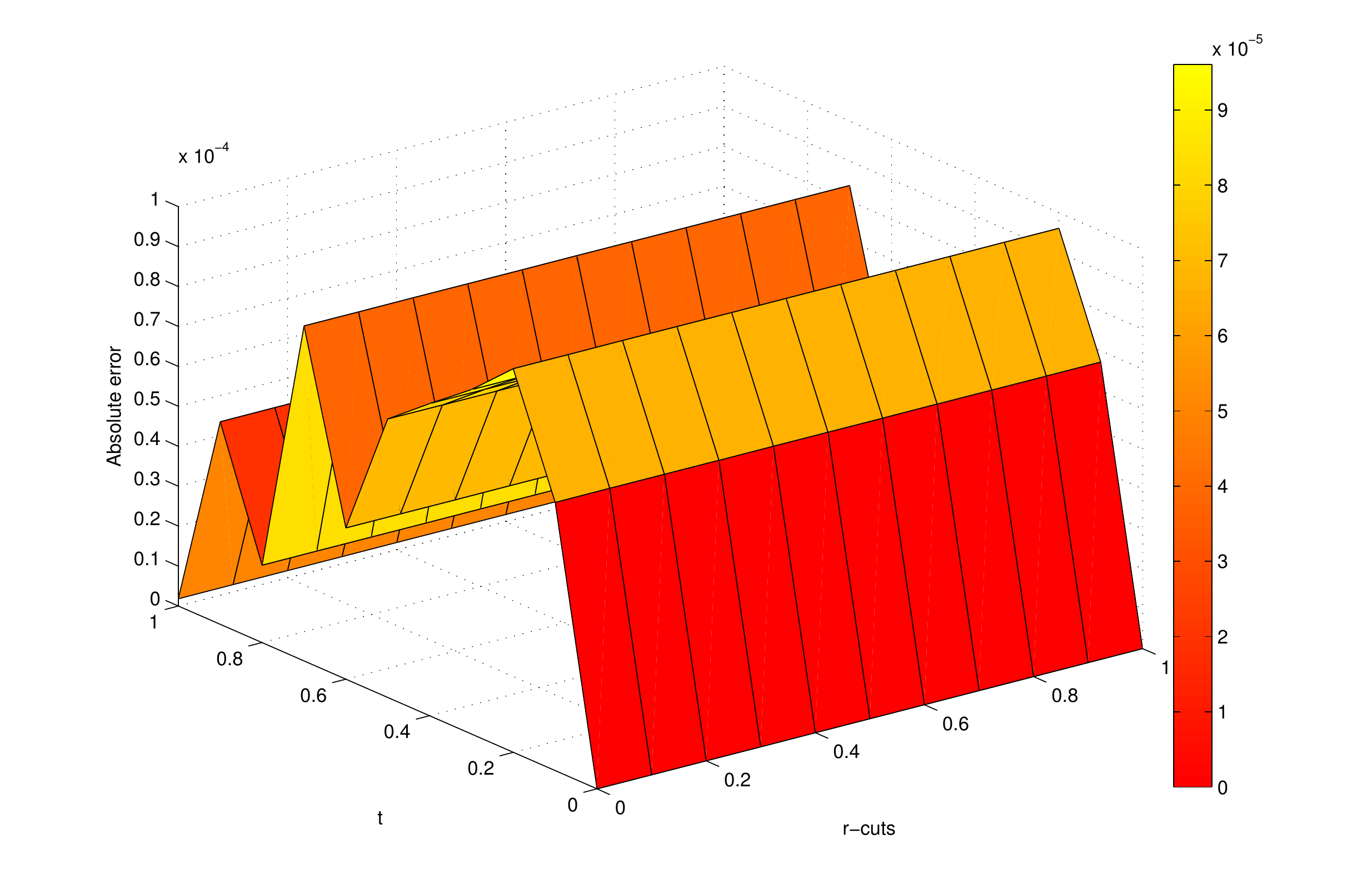}
}
\caption{\small{Absolute error functions (lower and upper bounds) of the problem 
\eqref{eq41} over $t\in[0,1]$, $r\in [0,1]$ 
and $v=\alpha=0.95$ with $N=9$ \cite{42}.}}
\label{fig111}
\end{figure}
\begin{figure}[!ht]
\subfigure[  $v=\alpha=0.78$
\label{subfig13a}]{%
\includegraphics[width=0.45\textwidth]{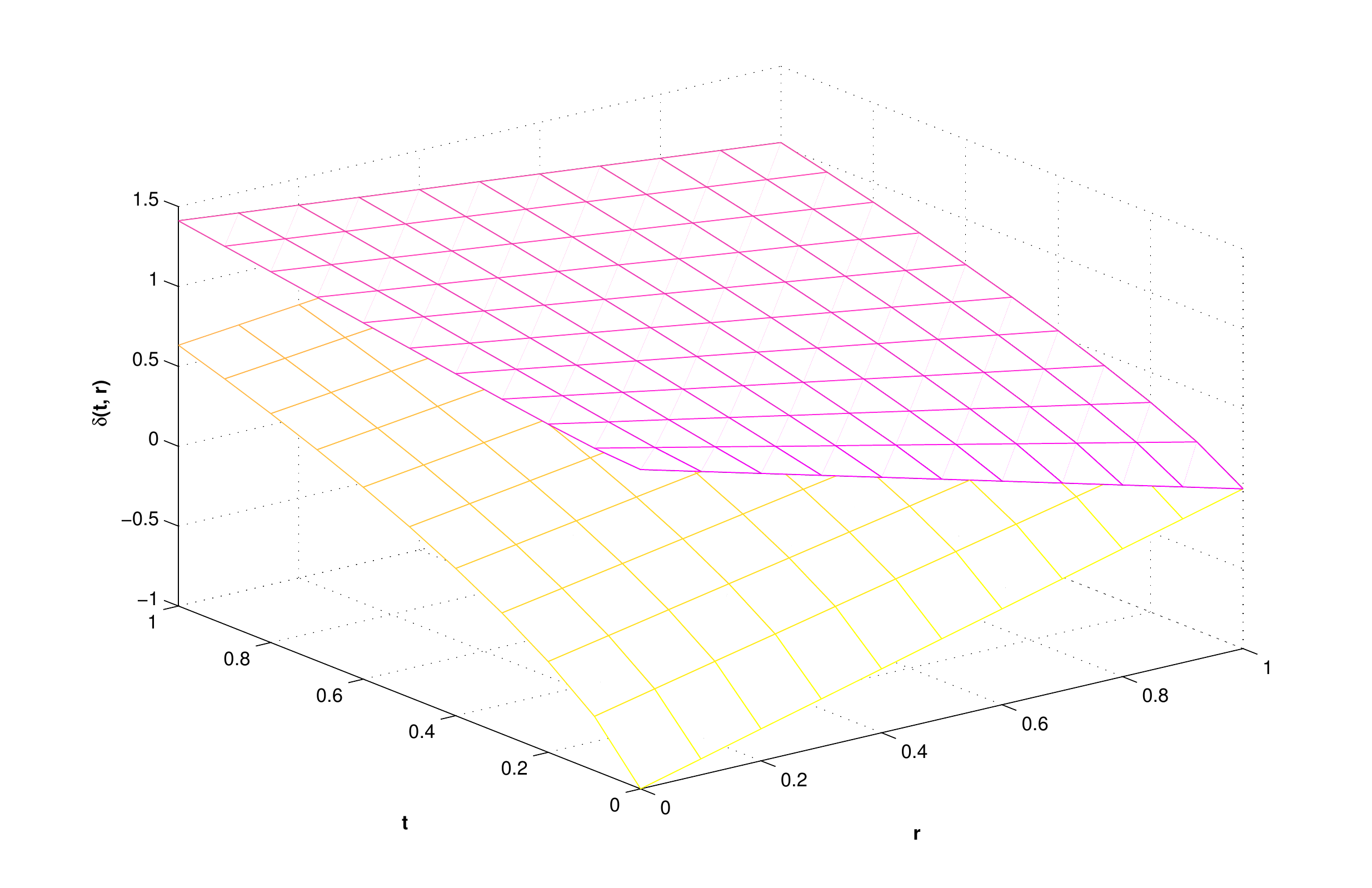}
}
\hfill
\subfigure[$v=\alpha=0.98$ \label{subfig13b}]{%
\includegraphics[ width=0.45\linewidth]{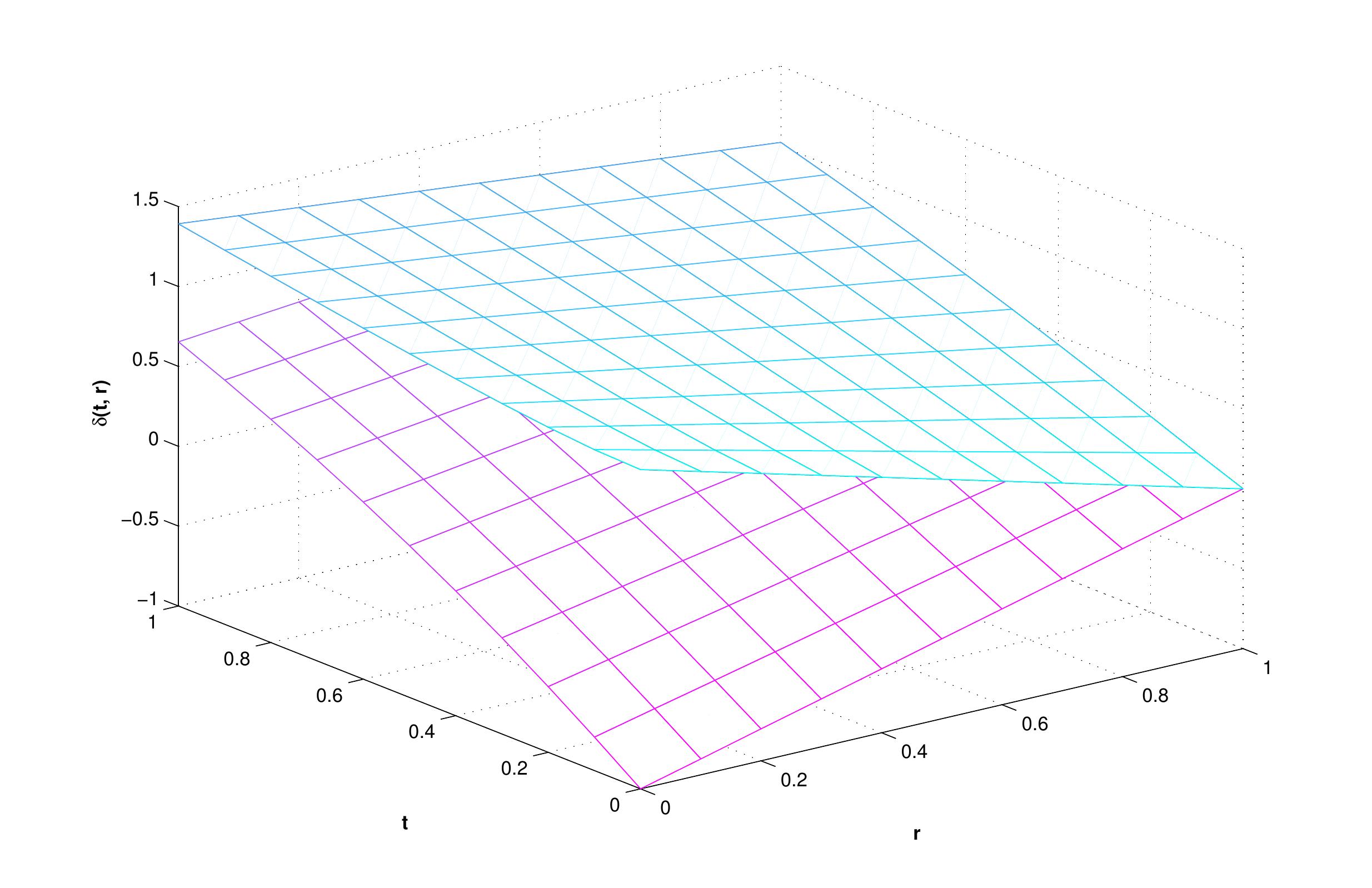}
}
\caption{\small{Fuzzy approximate solutions of the problem \eqref{eq41} 
for two different orders of the fractional derivative with $N=8$ 
over $t\in[0,1]$ and $r\in [0,1]$ \cite{42}.}}
\label{fig131}
\end{figure}


\section{Controllability}

Using a fixed point strategy is one of the fruitful methods to establish 
the controllability of nonlinear dynamical control system. In the last few 
years, some interesting and important controllability results, concerning 
semilinear differential systems involving Caputo fractional derivatives,
were proved. Consider the following Sobolev-type fractional evolution
system:
\begin{equation}
\label{sy-4.2}
\left\{
\begin{aligned}
&^{C}\!D^{q}(Ex(t))=Ax(t)+Ef(t,x(t))+EBu(t), && t\in J=[0,a],\\
&Ex(0)=Ex_0, && x_0\in D(E),
\end{aligned}
\right.
\end{equation}
where $^{C}\!D^{q}$ is the Caputo fractional derivative of
order $0<q<1$, $A: D(A)\subset X\rightarrow X$. Here $X$ is a separable
Banach space with the norm $|\cdot|$; $E:D(E)\subset X \rightarrow X$
are two closed linear operators and the pair $(A,E)$ generates an
exponentially bounded propagation family $\{W(t)\}_{t\geq 0}$ of $D(E)$ to $X$. 
The state $x(\cdot)$ takes values in $X$ and the control function 
$u(\cdot)$ is given in $\mathcal{U}$, the Banach space 
of admissible control functions, where
\[
\mathcal{U}=\left\{
\begin{aligned}
&L^{p}(J, U),&&~\mbox{for}~q\in \left(\frac{1}{p},1\right)~\mbox{with}~1<p<\infty,\\
&L^{\infty}(J, U),&&~\mbox{for}~q\in (0,1),
\end{aligned}
\right.
\]
and $U$ is a Banach space with the norm $|\cdot|_U$. 
Operator $B$ is bounded and linear from $U$ into $D(E)$, 
and $f : J\times X \rightarrow D(E)\subset X$ will
be specified later.


\subsection{Characteristic solution operators}

We recall the concept of exponentially bounded
propagation family.

\begin{definition}[See \cite{Liang}]
\label{def-W(t)}
A strongly continuous operator family $\{W(t)\}_{t\geq 0}$ of $D(E)$
to a Banach space $X$ such that $\{W(t)\}_{t\geq 0}$ is
exponentially bounded, which means that there exist $\omega>0$ and
$M>0$ such that $|W(t)x|\leq Me^{\omega t}|x|$ for any $x\in D(E)$
and $t\geq0$, is called an exponentially bounded propagation
family for the abstract degenerate Cauchy problem
\begin{equation}
\label{Deg-Cauchy}
\left\{
\begin{aligned}
&(Ex(t))'=Ax(t),&&t\in J,\\
&Ex(0)=Ex_0,&&x_0\in D(E),
\end{aligned}
\right.
\end{equation}
if
\begin{equation}
\label{Deg-Cauchy-Inte}
(\lambda E-A)^{-1}Ex=\int_{0}^{\infty}e^{-\lambda t}W(t)xdt,~x\in
D(E)
\end{equation}
when $\lambda>\omega$. In this case, we say that problem \eqref{Deg-Cauchy} has an
exponentially bounded propagation family $\{W(t)\}_{t\geq 0}$.
Moreover, if \eqref{Deg-Cauchy-Inte}  holds, then we also say that the
pair $(A,E)$ generates an exponentially bounded propagation family
$\{W(t)\}_{t\geq 0}$.
\end{definition} 

Let
\begin{equation}
\label{two-new-operators}
\begin{aligned}
\mathscr{T}_{(A,E)}(t)&=\int_{0}^{\infty}\Psi_{q}(\theta)W(t^{q}\theta)d\theta,\\
\mathscr{S}_{(A,E)}(t)&=q\int_{0}^{\infty}\theta
\Psi_{q}(\theta)W(t^{q}\theta)d\theta,
\end{aligned}
\end{equation}
where $\Psi_{q}(\theta)$ is the Wright function
\begin{equation*}
\Psi_\alpha(\theta)=\sum_{n=0}^{\infty}\frac{(-\theta)^n}{n!\Gamma(-\alpha n+1-\alpha)},
\quad \theta\in \mathbb{C}
\end{equation*}
with $0<\alpha<1$. We can introduce the following definition of mild solution for
system \eqref{sy-4.2}.

\begin{definition}
\label{new-def-sol}
For each $u\in \mathcal{U}$ and $x_0\in D(E)$, by a mild solution of
system \eqref{sy-4.2} we mean a function $x\in C(J,X)$ satisfying
\begin{align*}
x(t)=&\mathscr{T}_{(A,E)}(t)x_0+\int_{0}^{t}(t-s)^{q-1}\mathscr{S}_{(A,E)}(t-s)f\left(s,x(s)\right)ds\\
&+\int_{0}^{t}(t-s)^{q-1}\mathscr{S}_{(A,E)}(t-s)Bu(s)ds,
\quad t\in J.
\end{align*}
\end{definition}


\subsection{Controllability results}
\label{4.3.3}

In this subsection, we study the controllability of
system \eqref{sy-4.2} by utilizing the theory of propagation family.

\begin{definition}
System \eqref{sy-4.2} is said to be controllable on the interval
$J$ if for every $x_0 \in D(E)$ and every $x_{1}\in D(E)$ there
exists a control $u\in \mathcal{U}$ such that the mild solution
$x$ of system \eqref{sy-4.2} satisfies $x(a) = x_{1}$.
\end{definition}

We pose the following assumptions:
\begin{enumerate}
\item[$(H_1)$] the pair $(A,E)$ generates an exponentially
bounded propagation family $\{W(t)\}_{t\geq 0}$ of $D(E)$ to $X$;

\item[$(H_2)$] $\{W(t)\}_{t\geq 0}$ is a norm continuous family for
$t>0$ and $\|W(t)\|_{\mathscr{L}(X)}\leq M_1$ for $t\geq 0$;

\item[$(H_3)$] the control function $u(\cdot)$ takes from
$\mathcal{U}$, the Banach space of admissible control functions,
either $\mathcal{U}=L^{p}(J, U)$ for $q\in (\frac{1}{p},1)$ with
$1<p<\infty$ or $\mathcal{U}=L^{\infty}(J, U)$ for $q\in (0,1)$,
where $U$ is Banach space;

\item[$(H_4)$] $B:U\to D(E)$ is a bounded linear operator and the
linear operator $\mathbb{W} : \mathcal{U}\to X$ defined by
$$
\mathbb{W}u=\int_{0}^{a}(a-s)^{q-1}\mathscr{S}_{(A,E)}(a-s) Bu(s) ds
$$
has a bounded right inverse operator $\mathbb{W}^{-1}: X\to \mathcal{U}$.
\end{enumerate}
It is easy to see that $\mathbb{W}u\in X$ and $\mathbb{W}$ is well
defined due to the following fact:
\begin{align*}
|\mathbb{W}u|=&\bigg|\int_{0}^{a}(a-s)^{q-1}\mathscr{S}_{(A,E)}(a-s)
Bu(s) ds\bigg|\\
\leq& \frac{M_1\|B\|_{\mathscr{L}(U,X)}}{\Gamma(q)}\int_{0}^{a}(a-s)^{q-1} |u(s)|_{U} ds\\
\leq& \left\{
\begin{aligned}
&\frac{M_1\|B\|_{\mathscr{L}(U,X)}}{\Gamma(q)}\left(\frac{p-1}{qp-1}
a^{\frac{qp-1}{p-1}}\right)^{\frac{p-1}{p}}\|u\|_{L^pJ},\\
&~~~~~~~~~~~~~~~~~\mbox{if}~q
\in\left(\frac{1}{p},1\right),~u\in \mathcal{U}=L^p(J,U),~1<p<\infty,\\
&\frac{M_1\|B\|_{\mathscr{L}(U,X)}a^q}{\Gamma(q+1)}\|u\|_{L^\infty J},\\
&~~~~~~~~~~~~~~~~~\mbox{if}~q
\in (0,1),~u\in \mathcal{U}=L^\infty(J,U).
\end{aligned}
\right.
\end{align*}
Meanwhile,
\begin{equation}
\int_{0}^{t}(t-s)^{q-1} |u(s)|_{U}ds \leq K_q \|u\|_{L^pJ}
\end{equation}
for any $t\in J$, where
\begin{equation*}
K_q = \left\{
\begin{aligned}
&\left(\frac{p-1}{qp-1}a^{\frac{qp-1}{p-1}}\right)^{\frac{p-1}{p}}\|u\|_{L^pJ},\\
&\quad\quad\quad\quad\quad~\mbox{if}~q
\in (\frac{1}{p},1),~u\in \mathcal{U}=L^p(J,U),~1<p<\infty,\\
&\frac{a^q}{q}\|u\|_{L^\infty J},~\mbox{if}~q\in (0,1),~u
\in \mathcal{U}=L^\infty(J,U).
\end{aligned}
\right.
\end{equation*}
Next, we assume the hypothesis
\begin{enumerate}
\item [$(H_5)$] $f$ satisfies the following two conditions:
\begin{enumerate}
\item[(i)] for each $x\in X$ the function $f(\cdot, x) : J \rightarrow
D(E)\subset X$ is strongly measurable and, for each $t\in J$, the
function $f(t,\cdot):X\rightarrow D(E)\subset X$ is continuous;

\item[(ii)] for each $k>0$, there is a measurable function $g_{k}$ such that
$$
\sup_{|x|\leq k}|f(t,x)|\leq g_{k}(t),~\mbox{with}~\|g_k\|_\infty=\sup_{s\in J}g_k(s)<\infty
$$
and, for some $\gamma>0$, there exists sufficiently large $k_0$ such that
$$
\sup_{t\in J}\int_{0}^{t}(t-s)^{q-1}g_{k}(s)ds\le \gamma k,~~{\rm for~}k>k_0;
$$

\item[(iii)] there exists a positive constant $L>0$ such that
\[
\alpha(f(t,D))\leq L \alpha(D)
\]
for any bounded set $D\subset X$ and $t\in J$ a.e., 
where $\alpha$ is the measure of noncompactness.
\end{enumerate}
\end{enumerate}

Based on our assumptions, it is suitable 
to define the following control formula for an arbitrary function $x(\cdot)$:
\begin{equation}
\label{control-formula-4.2}
u(t)=\mathbb{W}^{-1}\left(x_{1}-\mathscr{T}_{(A,E)}(a)x_0
-\int_{0}^{a}(a-s)^{q-1}\mathscr{S}_{(A,E)}(a-s)f(s,x(s))ds\right).
\end{equation}

\begin{theorem}[See \cite{2014-2}]
\label{main-th1}
Assume $(H_1)$--$(H_5)$ are satisfied. Furthermore, assume that
\begin{equation}
\label{contraction-cond}
\rho=\left\{
\begin{aligned}
&\frac{\gamma M_{1}}{\Gamma(q)}\left(1
+\frac{a^{\frac{1}{2}}M_{1}\|B\|_{\mathscr{L}(U,X)}
K_q\|\mathbb{W}^{-1}\|_{\mathscr{L}(X,\mathcal{U})}}{\Gamma(q)}\right)<1,\\
&\qquad\qquad\qquad\qquad\qquad\qquad\qquad\qquad~\mbox{if~~}\mathcal{U}=L^2(J,U),\\
&\frac{\gamma M_{1}}{\Gamma(q)}\left(1+\frac{M_{1}\|B\|_{\mathscr{L}(U,X)}K_q\|
\mathbb{W}^{-1}\|_{\mathscr{L}(X,\mathcal{U})}}{\Gamma(q)}\right)<1,\\
&\qquad\qquad\qquad\qquad\qquad\qquad\qquad\qquad~\mbox{if~~}\mathcal{U}
=L^\infty(J,U),
\end{aligned}
\right.
\end{equation}
and
\begin{equation}
\label{MNC-condition}
\ell L\big(1+\ell\|B\|_{\mathscr{L}(U,X)}\|
\mathbb{W}^{-1}\|_{\mathscr{L}(X,\mathcal{U})}\big)<1,
\end{equation}
where $\ell=\frac{a^{q}M_1}{\Gamma(q+1)}$.
Then system \eqref{sy-4.2} is controllable on $J$.
\end{theorem}

\begin{corollary}[See \cite{2014-2}]
\label{main-cor1}
Let the assumptions in Theorem~\ref{main-th1} be satisfied. 
The set of mild solutions of system \eqref{sy-4.2} is a nonempty 
and compact subset of $C(J,X)$ with $u(t)$ given by \eqref{control-formula-4.2}.
\end{corollary}


\section{Approximate controllability}

Let $X$ be a Hilbert space with a scalar product
$\langle\cdot,\cdot\rangle$ and the corresponding norm $|\cdot|$.
We consider the following Sobolev-type fractional evolution system:
\begin{equation}
\label{sy.1}
\left\{
\begin{aligned}
&^{C}\!D^{q}(Ex(t))+Ax(t)=f(t,x(t))+Bu(t),
\quad t\in J=[0,a],\\
&x(0)+\sum^{m}_{k=1}a_k{x(t_k)}=0,
\end{aligned}
\right.
\end{equation}
where $^{C}\!D^{q}$ is the Caputo fractional
derivative of order $0<q<1$, $E$ and $A$ are two linear operators 
with domains contained in  $X$ and ranges still contained in  $X$, 
the pre-fixed points $t_k$ satisfy $0=t_0<{t_1}<{t_2}<\cdots<{t_m}<t_{m+1}=a$ 
and $a_k$ are real numbers.

In order to guarantee that $-AE^{-1}:X\to X$ generate a semigroup
$\{W(t)\}_{t\geq 0}$, we consider that the operators $A$ and $E$
satisfy the following conditions:
\begin{enumerate}
\item [$(S_1)$] $A: D(A)\subset X\rightarrow X$ 
and $E : D(E)\subset X \rightarrow X$ are linear, $A$ is closed;

\item [$(S_2)$] $D(E) \subset D(A)$ and $E$ is bijective;

\item [$(S_3)$] $E^{-1}: X\to D(E)$ is compact;

\item [$(S_3)'$] $E^{-1}: X\to D(E)$ is continuous.
\end{enumerate}
Now we note that
\begin{enumerate}
\item [(i)] $(S_3)'$ implies that $E$ is closed;

\item [(ii)] $(S_3)$ implies $(S_3)'$;

\item [(iii)] it follows from $(S_1)$, $(S_2)$, $(S_3)'$ and the closed graph theorem 
that $-AE^{-1}:X\to X$ is bounded, which generates a uniformly continuous
semigroup $\{W(t)\}_{t\geq 0}$ of bounded linear operators from $X$ to itself.
\end{enumerate}
Denote by $\rho(-AE^{-1})$ the resolvent set of $-AE^{-1}$. If we
assume that the resolvent $R(\lambda;-AE^{-1})$ is compact, then
$\{W(t)\}_{t\geq 0}$ is a compact semigroup (see \cite{3-pa}).

The state $x(t)$ takes values in $X$ and the control function
$u(\cdot)$ is given in $\mathcal{U}$, the Banach space of admissible control
functions, where $\mathcal{U}= L^{p}(J, U)$ for 
$q\in \left(\frac{1}{p},1\right)$ with $1<p<\infty$ and $U$ is a Hilbert space. 
Moreover, $B\in\mathscr{L}(U,X)$ is a bounded linear operator and $f :
J\times X \rightarrow X$ will be specified later.

Define the following two operators:
\begin{equation}
\label{two-new-operators1}
\begin{aligned}
&\mathscr{T}_{(A,E)}(t)=\int_{0}^{\infty}\Psi_{q}(\theta)W(t^{q}\theta)d\theta,\\
&\mathscr{S}_{(A,E)}(t)=q\int_{0}^{\infty}\theta
\Psi_{q}(\theta)W(t^{q}\theta)d\theta,
\end{aligned}
\end{equation}
where $\Psi_{q}(\theta)$ is the Wright function.
Now we introduce Green function:
\begin{equation}\label{green}
\begin{aligned}
G_{(A,E)}(t,s)&=E^{-1}G^0_{(A,E)}(t,s)\\
&=E^{-1}\left(-\displaystyle\sum^{m}_{k=1}\mathscr{T}_{(A,E)}(t)\chi_k(s) 
\Theta(t_k-s)^{q-1}\mathscr{S}_{(A,E)}(t_k-s)
+\chi_t(s)(t-s)^{q-1}\mathscr{S}_{(A,E)}(t-s)\right),
\end{aligned}
\end{equation}
for $t,s\in J$, where
$$
\chi_k(s)=\left\{\begin{aligned}&a_k,&&~{\rm for}~s\in[0,t_k),\\
&0, &&~{\rm for}~s\in[t_k,a],
\end{aligned}\right.
\qquad \chi_t(s)=\left\{\begin{aligned}&1,&&~{\rm for}~s\in[0,t),\\
&0,&&~{\rm for}~s\in[t,a].
\end{aligned}\right.
$$
Hence, we have that $\chi_k(s)(t_k-s)^{q-1}=0$ for $s\in[t_k,a]$ and
$\chi_t(s)(t-s)^{q-1}=0$ for $s\in[t,a]$.

Now, we introduce a suitable definition of mild solution.

\begin{definition}
\label{Definition 2.1ms} 
For each $u\in \mathcal{U}$,  by a mild solution of
system \eqref{sy.1} we mean a function $x\in C(J,X)$ satisfying
\begin{align*}
x(t)=\int_0^aG_{(A,E)}(t,s)\big(f(s,x(s))+Bu(s)\big)ds,
\quad t\in J.
\end{align*}
\end{definition}


\subsection{Linear systems}
\label{4.4.2}

Consider the following linear system:
\begin{equation}
\label{sy.1-linear}
\left\{
\begin{aligned}
&^{C}_{0}\!D^{q}_{t}(Ex(t))=Ax(t)+Bu(t),&&t\in J,\\
&x(0)+\sum^{m}_{k=1}a_k{x(t_k)}=0.
\end{aligned}
\right.
\end{equation}
Using the mild solution of \eqref{sy.1-linear}, we get
\begin{equation*}
x(a)=\int_0^aG_{(A,E)}(a,s)Bu(s)ds.
\end{equation*}
Now we recall the following result.

\begin{theorem}[See \cite{Mahmudov}]
\label{Theorem 2.1}  
Assume that $\Gamma: X\to X$ is symmetric. 
Then the following two conditions are equivalent:
\begin{enumerate}
\item[{\rm (i)}] $\Gamma: X\to X$ is positive, that is, $\langle x,\Gamma
x\rangle>0$ for all nonzero $x\in X$;

\item[{\rm (ii)}] for all $\eta\in X$, 
$x_{\varepsilon}(\eta)=\varepsilon(\varepsilon I+\Gamma)^{-1}(\eta)$ 
strongly converges to zero as $\varepsilon\to 0+$.
\end{enumerate}
\end{theorem}

We apply Theorem~\ref{Theorem 2.1} with $\Gamma^a_0$. Then, we have
\begin{align*}
\langle x^*,\Gamma^a_{0}x^*\rangle
&=\left\langle x^*,\int_{0}^{a}G_{(A,E)}(a,s)BB^*G^*_{(A,E)}(a,s)dsx^*\right\rangle\\
&=\int_{0}^{a}\left|B^*G^*_{(A,E)}(a,s)x^*\right|^2ds\\
&=\int_{0}^{a}\left|(P^*x^*)(s)\right|^2ds
\end{align*}
for any $x^*\in X$. Note that 
$P^* : X\to \mathcal{U}^*=L^{p^*}(J, U)\subset L^{p}(J, U)
=\mathcal{U}\subset L^{2}(J, U)$, since $1<p\le 2$. 
So the above last integral is well defined. 
We also get that $\langle x^*,\Gamma^a_{0}x^*\rangle>0$ if
and only if $P^*x^*\ne0$, i.e., $x^*\notin \ker P^*$. Consequently,
$\Gamma_0^a$ is positive if and only if $\ker P^*=\{0\}$, i.e.,
$\Gamma_0^a$ is positive if and only if the linear system
\eqref{sy.1-linear} is approximately controllable on $J$. Setting
$$
R(\varepsilon;\Gamma^a_{0})=(\varepsilon I+\Gamma^a_{0})^{-1}:X\to X, 
\quad \varepsilon>0,
$$
we arrive at the following result by Theorem~\ref{Theorem 2.1} 
(see also \cite{Mahmudov}).

\begin{theorem}[See \cite{jota-2015}]
\label{Theorem 2.2} 
Let $\frac{1}{2}<q<1$. The linear system \eqref{sy.1-linear} 
is approximately controllable on $J$ if and only if 
$\varepsilon R(\varepsilon;\Gamma^a_{0})\to 0$ as
$\varepsilon\to 0+$ in the strong topology.
\end{theorem}

Finally, we note that $R(\varepsilon;\Gamma_0^a)$ is continuous with
$\|R(\varepsilon;\Gamma_0^a)\|_{\mathscr{L}(X)}\le \frac{1}{\varepsilon}$.


\subsection{Approximate controllability}
\label{4.4.3}

In this subsection, we study the approximate controllability of
system \eqref{sy.1} by imposing that the corresponding linear system is
approximately controllable.

\begin{definition}
\label{Def. 3.1} 
Let $x(a;x(0),u)$ be the state value of
system \eqref{sy.1} at terminal time $a$ corresponding to the
control $u\in \mathcal{U}$ and nonlocal initial condition $x(0)$. 
System \eqref{sy.1} is said to be approximately controllable on the
interval $J$ if the closure $\overline{\mathfrak{R}(a,x(0))}=X$.
Here, $\mathfrak{R}(a,x(0))=\{x(a;x(0), u):u\in \mathcal{U}\}$ 
is called the reachability set of system \eqref{sy.1} at terminal time $a$.
\end{definition}

In the sequel, we introduce the following assumptions:
\begin{enumerate}
\item[$(H_1)$] $(S_1)$, $(S_2)$, and $(S_3)$ hold;

\item[$(H_2)$] $f: J\times X\to X$ is continuous such that
\[
g_k=\sup_{t\in J,|x|\leq
k}|f(t,x)|<\infty \quad \mbox{with}\quad 
\liminf_{k\to \infty}\frac{g_k}{k}=0;
\]

\item[$(H_3)$] system \eqref{sy.1-linear} is approximately controllable on $J$.
\end{enumerate}

Recalling condition $(H_3)$ and Theorem~\ref{Theorem 2.2}, 
we define the following control formula for any 
$x\in C(J,X)$ and $h\in X$:
\begin{equation}
\label{control-formula}
u_{\varepsilon}(t,x)=B^{*}G^*_{(A,E)}(a,t)
R(\varepsilon;\Gamma^a_{0})\Upsilon(x)
\end{equation}
with
\begin{equation*}
\Upsilon(x)=h-\int_0^aG_{(A,E)}(a,s)f(s,x(s))ds.
\end{equation*}
For each $k>0$, define
\[
\mathcal {B}_{k}=\{x\in C(J,X): ~\|x\|\leq k\}.
\]
Of course, $\mathcal {B}_{k}$ is a bounded, closed, convex subset in
$C(J,X)$, which is Banach space with the norm $\|\cdot\|$.

\begin{theorem}[See \cite{jota-2015}]
\label{Theorem 5.5.1} 
Let $\frac{1}{2}<q<1$. Under assumptions
$(H_1)$--$(H_3)$, for any $\varepsilon>0$, 
there exists a $k(\varepsilon)>0$ such that 
${P}_{\varepsilon}$ has a fixed point in $\mathcal{B}_{k(\varepsilon)}$.
\end{theorem}

Now we present a main result.

\begin{theorem}[See \cite{jota-2015}]
\label{Theorem 5.5.2} 
Let  all the assumptions in Theorem \ref{Theorem 5.5.1} be
satisfied. Moreover, there exists $r$ with $rq>1$ and 
$N\in L^{r}(J,\mathbb{R}^+)$ such that $|f (t, x)|\leq N(t)$ 
for all $(t, x) \in J \times X$. Then system \eqref{sy.1} 
is approximately controllable on the interval $J$.
\end{theorem}


\section{Existence and optimal control}
\label{sec:exist_oc}

Consider the following nonlinear fractional 
finite time delay evolution system:
\begin{align}
\label{3.1.1}
\left\{
\begin{aligned}
&{^C\!D^{q}}x(t)=Ax(t)+f(t,x_{t},x(t))+B(t)u(t),
&&0< t\leq T,\\
&x(t)=\varphi(t), && -r\leq t\leq 0,
\end{aligned}
\right.
\end{align}
where ${^C\!D^{q}}$ denotes the Caputo fractional derivative 
of order $q\in(0,1)$, $A$ is the generator of a $C_{0}$-semigroup
$\{T(t)\}_{t\geq 0}$ on a Banach space $X$, $f$ is a $X$-value function,
$u$ takes value from another Banach space $Y$, and $B$ is a linear operator
from $Y$ into $X$. Define $x_{t}$ by $z_t(\theta)=z(t+\theta)$,
$\theta\in[-r,0]$. Functions $f$, $x_{t}$ and $\varphi$ are given 
and satisfy some conditions that will be specified later.

Throughout this section, let $X$ and $Y$ be two Banach spaces,
with the norms $|\cdot|$ and $|\cdot|_Y$, respectively.
By $\mathscr{L}(X,Y)$ we denote the space of bounded linear operators from $X$
to $Y$ equipped with the norm $\|\cdot\|_{\mathscr{L}(X,Y)}$. In particular,
when $X=Y$, then $\mathscr{L}(X,Y)=\mathscr{L}(X,X)=\mathscr{L}(X)$ and
$\|\cdot\|_{\mathscr{L}(X,Y)}=\|\cdot\|_{\mathscr{L}(X,X)}=\|\cdot\|_{\mathscr{L}(X)}$.
Suppose $r>0$ and $T>0$. Denote $J=[0,T]$ and
$M=\sup_{t\in J}\|T(t)\|_{\mathscr{L}(X)}$, which is a finite number. 
Let $C([-r,a],X)$, $a\geq 0$, be the Banach space of
continuous functions from $[-r,a]$ to $X$ with the usual sup-norm. 
For brevity, we denote $C([-r,a],X)$ simply by $C_{-r,a}$
and its norm by $\|\cdot\|_{-r,a}$. If $a=T$, then we denote this space
by $C_{-r,T}$ and its norm by $\|\cdot\|_{-r,T}$. If $a=0$, then we
denote this space by $C_{-r,0}$ and its norm by
$\|\cdot\|_{-r,0}$. For any $x\in C_{-r,T}$ and $t\in J$, 
define $x_{t}(s)=x(t+s)$ for $-r\leq s \leq 0$. Then 
$x_{t}\in C_{-r,0}$.


\subsection{Existence and uniqueness}
\label{4.2.2}

We make the following assumptions:
\begin{enumerate}
\item[$(H_1)$] $f:J\times C_{-r,0} \times X\rightarrow X$ satisfies:
\begin{enumerate}
\item [(i)] for each $x_{t}\in C_{-r,0}$, $x\in X$, $t\rightarrow
f(t,x_{t},x(t))$ is measurable;

\item [(ii)] for arbitrary $\xi_{1}$, $\xi_{2}\in C_{-r,0}$, $\eta_{1}$,
$\eta_{2}\in X$ satisfying $\|\xi_{1}\|_{-r,0}$,
$\|\xi_{2}\|_{-r,0}$, $|\eta_{1}|$, $|\eta_{2}|\leq\rho$, there
exists a constant $L_{f}(\rho)>0$ such that
\[
|f(t,\xi_{1},\eta_{1})-f(t,\xi_{2},\eta_{2})|\leq
L_{f}(\rho)(\|\xi_{1}-\xi_{2}\|_{-r,0}+|\eta_{1}-\eta_{2}|)
\]
for all $t\in J$;

\item [(iii)] there exists a constant $a_{f}>0$ such that
\[
|f(t,\xi,\eta)|\leq a_{f}(1+\|\xi\|_{-r,0}+|\eta|),
\mbox{ for all } \xi\in C_{-r,0},\ \eta\in X, \ t\in J.
\]
\end{enumerate}

\item[$(H_2)$] Let $Y$ be a reflexive Banach space from which
the controls $u$ take their values. The operator 
$B\in L_{\infty}(J,\mathscr{L}(Y,X))$, where 
$\|B\|_{\infty}$ stands for the norm of operator $B$ 
on the Banach space $L_{\infty}(J,\mathscr{L}(Y,X))$.

\item[$(H_3)$] The multivalued map $U\left( \cdot \right): 
J \rightarrow P(Y)$ has closed, convex, and bounded values
with $U\left(\cdot \right)$ graph measurable and 
$U\left( \cdot \right) \subseteq \Omega$, 
where $\Omega$ is a bounded set in $Y$.
\end{enumerate}

Introduce the admissible set
\[
U_{ad}=\left\{ v\left( \cdot \right): J\rightarrow 
Y\mbox{ strongly measurable}, v(t)\in U(t) \mbox{ a.e.}~t\in J\right\}.
\]
Obviously, $U_{ad}\neq \emptyset$ (see Theorem~2.1 of \cite{Zeidler-1990}) 
and $U_{ad}$ $\subset L^{p}(J,Y)(1< p<+\infty)$ is bounded, closed, and convex. 
It is obvious that $Bu\in L^{p}(J,X)$ for all $u\in U_{ad}$.

We give the following definition of mild
solution for the problem below.

\begin{definition} 
\label{Definition 3.1.}
For any $u\in L^{p}(J,Y)$, if there exist $T=T(u)>0$ 
and $x\in C([-r,T],X)$ such that
\begin{equation}
\label{equ-3.2.1}
x(t)=\left\{
\begin{aligned}
&S_{q}(t)\varphi(0)+\int_{0}^{t}(t-s)^{q-1}P_{q}(t-s)f\left(s,x_{s},x(s)\right)ds\\
&+\int_{0}^{t}(t-s)^{q-1}P_{q}(t-s)B(s)u(s)ds,\quad~0\leq t\leq T,\\
&\varphi(t), \quad\quad\quad\quad\quad\quad\quad\quad\quad\quad\quad\quad\quad\quad\quad-r\leq t\leq 0,
\end{aligned}\right.
\end{equation}
then system \eqref{3.1.1} is called mildly solvable with respect
to $u$ on $[-r,T]$, where
\[
S_{q}(t)=\int_{0}^{\infty}\Psi_{q}(\theta)T(t^{q}\theta)d\theta,\quad
P_{q}(t)=q\int_{0}^{\infty}\theta\Psi_{q}(\theta)T(t^{q}\theta)d\theta,
\]
and $\Psi_{q}(\theta)$ is the Wright function.
\end{definition}

\begin{theorem}[See \cite{new-0-32}]
\label{Theorem 3.1.}
Assume that $(H_1)$, $(H_2)$, and $(H_3)$ hold. Then for each $u\in U_{ad}$
and for some $p$ such that $pq>1$, system \eqref{3.1.1} is mildly
solvable on $[-r,T]$ with respect to $u$ and the mild solution is unique.
\end{theorem}


\subsection{Optimal control} 
\label{4.2.4}

In what follows, we consider the fractional optimal
control of system \eqref{3.1.1}. Precisely, we consider 
the following optimal control problem in Lagrange form:
find a control $u^{0} \in U_{ad}$ such that
\begin{equation}
\label{P}
\tag{$P$}
\begin{gathered}
J(u^{0}) \leq J(u) \mbox{ for all } u\in U_{ad},\\
J(u) =\int_{0}^{T}\mathcal {L}( t, x^{u}_{t}, x^{u}(t), u(t)) dt.
\end{gathered}
\end{equation}
Here $x^{u}$ denotes the mild solution of system \eqref{3.1.1}
corresponding to the control $u\in U_{ad}$.
	
For the existence of solution to problem \eqref{P}, we introduce the
following assumption:
\begin{enumerate}
\item[$(H_4)$]
\begin{enumerate}
\item [(i)] functional
$\mathcal {L}: J\times C_{-r,0}\times
X\times Y\to \mathbb{R}\cup \{ \infty \}$ 
is Borel measurable;

\item [(ii)] $\mathcal {L}(t,\cdot,\cdot,\cdot )$ is sequentially lower
semicontinuous on $ C_{-r,0}\times X \times Y$ for almost all $t\in J$;

\item [(iii)] $\mathcal {L}(t,x,y,\cdot )$ is convex on $Y$ for each 
$x\in C_{-r,0}$, $y\in X$ and almost all $t\in J$;

\item [(iv)] there exist constants $d, e\geq 0$, $j>0$, 
$\varphi$ is nonnegative, and $\varphi \in L^{1}( J,\mathbb{R})$ such that
\[
\mathcal {L}(t,x,y,u)\geq \varphi (t)+d\| x\|_{-r,0} +e|y| 
+j\|u\|^{p}_{Y}.
\]
\end{enumerate}
\end{enumerate}

Now we can give the following result on existence of fractional
optimal controls for problem \eqref{P}.

\begin{theorem}[See \cite{new-0-32}]
\label{Theorem 3.3.} 
Under the assumptions of
Theorem \ref{Theorem 3.1.} and $(H_4)$, suppose that $B$ is a
strongly continuous operator. Then the optimal control problem \eqref{P}
admits at least one optimal pair, i.e., there exists an admissible
control $u^{0} \in U_{ad}$ such that
\[
J(u^{0})=\int_{0}^{T}\mathcal {L}( t, x^{0}_{t}, x^{0}(t),
u^{0}(t)) dt \leq J(u) \quad \mbox{  for } u \in U_{ad}.
\]
\end{theorem}

\begin{remark}
\label{Remark 3.4.}
Condition $(H_3)$ in Theorem~\ref{Theorem 3.3.} can be replaced
by the following condition:
\begin{enumerate}
\item [$(H_3)'$] $\mathscr{U}$ is a weakly compact subset of $Y$ and
$t\to U(t)$ is a map with measurable values in $P_{cl,cv}(\mathscr{U})$.
\end{enumerate}
\end{remark}

\begin{theorem}[See \cite{new-0-32}]
\label{Theorem 3.4.}
Under the assumptions in Theorem~\ref{Theorem 3.3.} with $(H_3)$
replaced by $(H_3)'$, let
$$
U_{ad}=\left\{u\left( \cdot \right) : J\rightarrow Y
\mbox{~is~ strongly~ measurable}, u(t)\in U(t), t\in J\right\}.
$$
Then there exists an optimal control for problem \eqref{P}.
\end{theorem}


\section{Optimal feedback control}
\label{sec7}

Control systems are often based on the feedback principle,
whereby the signal to be controlled is compared with a desired reference
signal and the discrepancy used to compute a corrective control action.
Consider the following semilinear fractional feedback control system:
\begin{align}
\label{S.1.1}
\left\{
\begin{aligned}
&^{C}\!D^{q}x(t)=Ax(t)+f(t,x(t),u(t)),&&t\in J=[0,T], \\
& x(0)=x_{0},
\end{aligned}
\right.
\end{align}
where $^{C}\!D^{q}$ is the Caputo fractional derivative of order
$q\in(0,1)$, and $A: D(A)\to X$ is the infinitesimal generator of
a compact $C_{0}$-semigroup $\{T(t)\}_{t\geq 0}$ in a
reflexive Banach space $X$. The control $u$ takes value from
$U[0,T]$, which is the control set, $f:J\times X\times U\to X$ will
be specified later.


\subsection{Existence of feasible pairs}

Denote by $X$ a reflexive Banach space with norm $|\cdot|$, and by
$U$ a Polish space which is a separable completely metrizable
topological space. Let $C(J,X)$ be the Banach space of continuous
functions from $J$ to $X$ with the usual sup-norm. Suppose that
$A: D(A)\to X$ is the infinitesimal generator of a compact
$C_{0}$-semigroup $\{T(t)\}_{t\geq 0}$. This means that
there exists $M>0$ such that
$\sup_{t\in J}\|T(t)\|_{\mathscr{L}(X)}\leq M$. By
\[
O_{r}(x)=\{y\in X:|y-x|\leq r\},
\]
we denote the ball centered at $x$ with the radius $r>0$.

\begin{definition}[See \cite{LiXunjing1}]
\label{def0}
Let $E$ and $F$ be two metric spaces. A multifunction
$\digamma:E\to P(F)$ is said to be pseudo-continuous at $t\in E$ if
\[
\bigcap_{\epsilon>0}\overline{\digamma(O_{\epsilon}(t))}=\digamma(t).
\]
We say that $\digamma$ is pseudo-continuous on $E$ if it is
pseudo-continuous at each point $t\in E$.
\end{definition}

We make the following assumptions:
\begin{enumerate}
\item [$(H_S)$] $X$ is a reflexive Banach space and $U$ is a Polish space;

\item [$(H_A)$] $A$ is the infinitesimal generator of a compact
$C_{0}$-semigroup $\{T(t)\}_{t\geq 0}$ on $X$;

\item [$(H_1)$] $f:J\times X\times U\rightarrow X$ is Borel measurable in
$(t,x,u)$ and is continuous in $(x,u)$;

\item [$(H_2)$] $f$ is local Lipschitz continuous with respect to $x$,
i.e., for any constant $\rho>0$, there is a constant $L(\rho)>0$
such that
\[
|f\left(t,x_{1},u\right) - f\left(t,x_{2},u\right) |
\leq L\left( \rho \right)| x_{1}-x_{2}|
\]
for every $x_{1}$, $x_{2}\in X$, $t\in J$, and uniformly with respect to 
$u\in U$ provided $|x_{1}|$, $|x_{2}|\leq\rho$;

\item [$(H_3)$] for arbitrary $t\in J$, $x\in X$, and $u\in U$, there exists a
positive constant $M>0$ such that
\[
| f\left( t,x,u\right) | \leq M(1+|x|);
\]

\item [$(H_4)$] for almost all $t\in J$, 
the set $f(t,x,\digamma(t,x))$ satisfies 
\[
\bigcap_{\delta>0}\overline{\rm co}f(t,O_{\delta}(x),\digamma(O_{\delta}(t,x)))
=f(t,x,\digamma(t,x));
\]

\item [$(H_U)$] $\digamma: J\times X\to P(U)$ is pseudo-continuous.
\end{enumerate}

Let $U[0,T]=\{u: J\to U, u(\cdot) \mbox{ is measurable}\}$.
Then any element in the set $U[0,T]$ is called a control on $J$.
In the following, we introduce the definition of mild
solution for system \eqref{S.1.1}.

\begin{definition}
A mild solution $x \in C(J,X)$ of system \eqref{S.1.1} is defined as
a solution of the integral equation
\begin{equation}
\label{sol-express}
x(t)=\mathscr{T}(t)x_{0}+\int_{0}^{t}(t-\theta)^{q-1}\mathscr{S}\left(t-\theta
\right)f\big(\theta ,x\left(\theta
\right),u(\theta)\big)d\theta,
\quad t\in J,
\end{equation}
where
\[
\mathscr{T}(t)=\int_{0}^{\infty}\Psi_{q}(\theta)T(t^{q}\theta)d\theta,
\quad \mathscr{S}(t)=q\int_{0}^{\infty}\theta\Psi_{q}(\theta)T(t^{q}\theta)d\theta.
\]
Here, $\Psi_{q}(\theta)$ is the Wright function.
\end{definition}

Any solution $x(\cdot)\in C(J,X)$ of system \eqref{S.1.1} is
referred as a state trajectory of the fractional evolution
equation corresponding to the initial state $x_{0}$ and 
the control $u(\cdot)$.

\begin{theorem}[See \cite{Wang-SCL}]
\label{Existence}
Assume $(H_S)$, $(H_A)$, $(H_1)$, $(H_2)$, and $(H_3)$ hold. 
Then there is a unique mild solution $x\in C(J,X)$ 
of system \eqref{S.1.1} for any $x_{0}\in X$ and $u\in U$, 
and
\[
\|x\|\leq M,
\]
for some constant $M>0$.
\end{theorem}

Next, we introduce the definition of feasible pair.

\begin{definition}
\label{D.2.1}
A pair $(x,u)$ is said to be feasible if $x$
satisfies \eqref{sol-express} and
\[
u(t)\in \digamma(t,x(t)),\quad a.e.~t\in J.
\]
\end{definition}

Let $[s,v]\subseteq J$,
\begin{align*}
\mathcal{H}[s,v]=&\{(x,u)\in C([s,v],X)\times U[s,v]:
(x,u)\text{~is~feasible}\},\\
\mathcal{H}[0,T]=&\{(x,u)\in C([0,T],X)\times U[0,T]: (x,u)\text{~is~feasible}\}.
\end{align*}
To solve the optimal feedback control problem, we need the following
result, which is an extension of the results corresponding to first
order semilinear evolution equations.

\begin{theorem}[See \cite{Wang-SCL}]
\label{Fe-Pair}
Assume that $(H_S)$, $(H_A)$, $(H_1)$--$(H_4)$ and $(H_U)$ hold. 
Then, for any $x_0\in X$ and $\frac{1}{p}<q<1$ and for some $p>1$, 
the set $\mathcal{H}[0,T]$ is nonempty:
\[
\mathcal{H}[0,T]\neq \emptyset.
\]
\end{theorem}


\subsection{Existence of optimal feedback control pairs}

We now consider the following Lagrange problem:
find a pair $(x^{0},u^{0}) \in \mathcal{H}[0,T]$ such that
\begin{equation}
\label{Q}
\tag{$Q$}
\begin{gathered}
\mathcal {J}(x^{0},u^{0}) \leq \mathcal {J}(x,u) 
\mbox{ for all }(x,u)\in \mathcal{H}[0,T],\\
\mathcal {J}(x,u)=\int_{0}^{T}\mathcal {L}(t\mbox{, }x(t)\mbox{, } u(t)) dt.
\end{gathered}
\end{equation}
We impose some assumptions on $\mathcal {L}$:
\begin{enumerate}
\item[$(L_1)$] functional $\mathcal {L}: J\times X\times U\to
R\cup \{ \infty \} $ is Borel measurable in $(t,x,u)$;

\item[$(L_2)$]  $\mathcal {L}(t,\cdot ,\cdot )$ is sequentially 
l.s.c. on $X\times U$ for almost all $t\in J$ and there is a
constant $M_{1}>0$ such that
\[
\mathcal {L}(t,x,u)\geq -M_{1},~(t,x,u)\in J\times X\times U.
\]
\end{enumerate}

For any $(t,x)\in J\times X$, we set
\begin{align*}
\mathcal{W}(t,x)=\{(z^{0},z)\in R\times X:~& z^{0}\geq \mathcal{L}(t,x,u),\\
&z=f(t,x,u),~u\in \digamma (t,x)\}.
\end{align*}
In order to prove existence of optimal control pairs for problem
\eqref{Q}, we assume that:
\begin{enumerate}
\item [$(H_C)$] for almost all $t\in J$, the map $\mathcal{W}(t,\cdot): X\to
P(R\times X)$ has the Cesari property, i.e.,
\[
\bigcap_{\delta>0}\overline{\rm co}\mathcal{W}(t,O_{\delta}(x))=\mathcal{W}(t,x)
\]
for all $x\in X$.
\end{enumerate}

\begin{theorem}[See \cite{Wang-SCL}]
\label{Exis-Pair}
Assume that the hypotheses $(H_S)$, $(H_A)$, $(H_1)$--$(H_4)$, $(H_U)$, $(L_1)$, 
$(L_2)$, and $(H_C)$ hold. Then the Lagrange problem \eqref{Q} 
admits at least one optimal control pair.
\end{theorem}

In \cite{MR2824730}, nonlocal problems for a class of fractional 
integrodifferential equations via fractional operators and optimal control 
in Banach spaces are investigated. The results make used of fractional calculus, 
H\"{o}lder inequality, $p$-mean continuity and fixed point theorems. 
Some existence results of mild solutions are also obtained \cite{MR2824730}.
Existence and uniqueness of solutions were proved in \cite{MR2886699},
by means of the H\"{o}lder inequality, a suitable singular Gronwall inequality, 
and fixed point theorems. See also \cite{MR3074465,Zhou-1,Zhou-2}.
Here we proceed by reviewing the main results of \cite{MyID:323}.


\section{Optimal solutions to relaxation in multiple control problems of Sobolev type}
\label{sec8}

We now address the optimality question to the relaxation in multiple control
problems described by Sobolev type nonlinear fractional differential equations
with nonlocal control conditions in Banach spaces. Moreover, we consider
the minimization problem of multi-integral functionals, with integrands
that are not convex in the controls, of control systems
with mixed nonconvex constraints on the controls. We prove,
under appropriate conditions, that the relaxation problem
admits optimal solutions. Furthermore, we show that those optimal
solutions are in fact limits of minimizing sequences of systems
with respect to the trajectory, multi-controls,
and the functional in suitable topologies.


\subsection{Optimal control problems with nonlocal nonlinear fractional differential equations}

Consider the following nonlocal nonlinear fractional control system of Sobolev type:
\begin{equation}
\label{eq:1.1}
L~^CD_t^{\alpha}[Mx(t)]+Ex(t)=f(t, x(t), B_{1}(t)u_{1}(t),\ldots, B_{r-1}(t)u_{r-1}(t)),~ t\in I
\end{equation}
\begin{equation}
\label{eq:1.2}
x(0)+h(x(t), B_{r}(t)u_{r}(t))=x_{0},
\end{equation}
with mixed nonconvex constraints on the controls
\begin{equation}
\label{eq:1.3}
u_{1}(t),\ldots, u_{r}(t)\in U(t, x(t))~ \text{a.e. on}~ I,
\end{equation}
where $^CD^{\alpha}_{t}$ is the Caputo fractional derivative of order
$\alpha$, $0<\alpha\leq1$, and $t\in I:=[0, a]$. Let $X, Y$ and $Z$
be three Banach spaces such that $Z$ is densely and continuously embedded
in $X$, the unknown function $x(\cdot)$ takes its values in $X$
and $x_{0}\in X$. We assume that the operators $E: D(E)\subset X\rightarrow Y$,
$M: D(M)\subset X\rightarrow Z$, $L: D(L)\subset Z\rightarrow Y$,
and $B_{1},\ldots,B_{r}:I\to\mathcal{L}(T, X)$ are linear and bounded
from $T$ into $X$. The space $T$ is a separable reflexive Banach space
modeling the control space. It is also assumed that $f: I\times X^{r}\rightarrow Y$
and $h: C(X^{2}, X)\rightarrow X$ are given abstract functions, to be specified later,
and $U:I\times X \rightrightarrows 2^T\backslash\{\emptyset\}$ is a multivalued map
with closed values, not necessarily convex.
Let $\widehat{\mathbb{R}} := \ ]-\infty, +\infty]$. For functions
$g_{1},\ldots,g_{r}:I\times X\times T \to \mathbb{R}$, we consider the problem
$$
\max\left\{J_{1},\ldots,J_{r}\Biggl\vert
\begin{array}{lll}
J_{1}(x, u_{1}):=\int_{I}g_{1}(t, x(t), u_{1}(t))dt\\
\qquad \vdots\\
J_{r}(x, u_{r}):=\int_{I}g_{r}(t, x(t), u_{r}(t))dt
\end{array}
\right\}
\longrightarrow \inf \eqno(P)
$$
on solutions of the control system \eqref{eq:1.1}--\eqref{eq:1.2} with constraint \eqref{eq:1.3}.
Let $g_{1, U},\ldots,g_{r, U}: I\times X\times T\to \widehat{\mathbb{R}}$ be the functions defined by
$$
\left\{
\begin{array}{ll}
g_{1,U}(t, x, u_1)
:=\left\{
\begin{array}{ll}
g_{1}(t, x, u_{1}),& \mbox{$u_{1}\in U(t, x)$},\\
+\infty, & \mbox{$u_{1}\notin U(t, x)$},
\end{array}\right.\\
\qquad \vdots\\
g_{r,U}(t, x, u_{r})
:=\left\{
\begin{array}{ll}
g_{r}(t, x, u_{r}),& \mbox{$u_{r}\in U(t, x)$},\\
+\infty, & \mbox{$u_{r}\notin U(t, x)$},
\end{array}\right.
\end{array}\right.
$$
and $g^{**}_{1}(t, x, u_{1}),\ldots,g^{**}_{r}(t, x, u_{r})$
be the bipolar of $u_{1}\to g_{1,U}(t, x, u_{1}),\ldots$,
$u_{r}\to g_{r,U}(t, x, u_{r})$, respectively. Along with problem
$(P)$, we also consider the relaxation problem
$$
\max\left\{J^{**}_{1},\ldots,J^{**}_{r}\Biggl\vert
\begin{array}{lll}
J^{**}_{1}(x, u_{1})=\int_{I}g^{**}_{1}(t, x(t), u_{1}(t))dt\\
\vdots\\
J^{**}_{r}(x, u_{r})=\int_{I}g^{**}_{r}(t, x(t), u_{r}(t))dt
\end{array}
\right\}
\longrightarrow \inf \eqno(RP)
$$
on the solutions of control system \eqref{eq:1.1}--\eqref{eq:1.2}
with the convexified constraints
\begin{equation}
\label{eq:1.4}
u_{1}(t),\ldots, u_{r}(t)\in \operatorname{cl} \operatorname{conv}U(t,x(t))\quad\text{a.e. on }I
\end{equation}
on the controls, where $conv$ denote the convex hull and $cl$ the closure.
In our results, we will denote by $\mathcal{R}_U$ and $\mathcal{T}r_{U},
(\mathcal{R}_{\operatorname{cl} \operatorname{conv}U}$
and $\mathcal{T}r_{\operatorname{cl} \operatorname{conv} U})$ the sets of all solutions and
all trajectories of control system \eqref{eq:1.1}--\eqref{eq:1.3}
(control system \eqref{eq:1.1}--\eqref{eq:1.2},\eqref{eq:1.4}, respectively).

We make the following assumptions:
\begin{itemize}
\item [(H$_1$)] $L: D(L)\subset Z\rightarrow Y$ and $M: D(M)\subset X\rightarrow Z$
are linear operators, and $E: D(E)\subset X\rightarrow Y$ is closed.
	
\item [(H$_2$)] $D(M)\subset D(E)$, $\text{Im}(M)\subset D(L)$ and $L$ and $M$ are bijective.
	
\item [(H$_3$)] $L^{-1}: Y\rightarrow D(L)\subset Z$
and $M^{-1}: Z\rightarrow D(M)\subset X$ are linear, bounded and compact operators.
\end{itemize}
Note that (H$_3$) implies that $L$ and $M$ are closed.
Indeed, if $L^{-1}$ and $M^{-1}$ are closed and injective, then
their inverse are also closed. From (H$_1$)--(H$_3$) and the closed graph theorem,
we obtain the boundedness of the linear operator $L^{-1}EM^{-1}: Z\rightarrow Z$. Consequently,
$L^{-1}EM^{-1}$ generates a semigroup $\lbrace Q(t), t\geq0\rbrace$, $Q(t):=e^{L^{-1}EM^{-1}t}$.
We assume that $M_{0}:=\sup_{t\geq0}\Vert Q(t)\Vert<\infty$ and, for short,
we denote $C_{1}:=\Vert L^{-1}\Vert$ and $C_{2}:=\Vert M^{-1}\Vert$.
According to previous definitions, it is suitable to rewrite problem
\eqref{eq:1.1}--\eqref{eq:1.2} as the equivalent integral equation
\begin{equation}
\label{eq:2.1}
Mx(t)=Mx(0)
+\frac{1}{\Gamma(\alpha)}\int_{0}^{t}(t-s)^{\alpha-1}
[-L^{-1}Ex(s)+L^{-1}f(s, x(s), B_{1}(s)u_{1}(s),\ldots, B_{r-1}(s)u_{r-1}(s))]ds,
\end{equation}
provided the integral in \eqref{eq:2.1} exists a.e. in $t\in J$.
Before formulating the definition of mild solution of system
\eqref{eq:1.1}--\eqref{eq:1.3}, we first introduce some necessary notions.
Let $I:=[0,a]$ be a closed interval of the real line with the Lebesgue
measure $\mu$ and the $\sigma$-algebra $\Sigma$ of $\mu$ measurable sets.
The norm of the space $X$ (or $T$) will be denoted by $\|\cdot\|_X$
(or $\|\cdot\|_T$). We denote by $C(I,X)$ the space of all continuous
functions from $I$ into $X$ with the supnorm given by
$\|x\|_{C}:=\sup_{t\in I}\|x(t)\|_X$ for $x\in C(I,X)$.
For any Banach space $V$, the symbol $\omega$-$V$ stands for $V$ equipped
with the weak topology $\sigma(V,V^*)$. The same notation will be used
for subsets of $V$. In all other cases, we assume that $V$ and its subsets
are equipped with the strong (normed) topology.

In what follows, $A:=-L^{-1}EM^{-1}: D(A)\subset Z\rightarrow Z$
is the infinitesimal generator of a compact analytic
semigroup of uniformly bounded linear operators $Q(\cdot)$ in $X$.
Then, there exists a constant $M_{0}\geq1$ such that $\Vert Q(t)\Vert\leq M_{0}$
for $t\geq0$. The operators $B_{i}\in L^{\infty}(I, \mathcal{L}(T, X))$,
and we let $\Vert B_{i}\Vert$ stand for $\Vert B_{i}\Vert_{L^{\infty}(I, \mathcal{L}(T, X))}$.

We make use of the following assumptions on the data of our problems.
\begin{itemize}
\item[(H1)] The nonlinear function $f: I\times X^{r}\to Y$ satisfies the following:
\begin{itemize}
\item[(1)] $t\to f(t, x_{1},\ldots,x_{r})$ is measurable
for all $(x_{1},\dots,x_{r})\in X^{r}$;
		
\item[(2)] $\Vert f(t, x_{1},\ldots,x_{r})-f(t, y_{1},\ldots,y_{r})\Vert_{Y}\leq k_{1}(t)
\sum_{i=1}^{r}\Vert x_{i}-y_{i}\Vert_{X}$ a.e. on $I$, $k_{1}\in L^{\infty}(I,\mathbb{R}^+)$;
		
\item[(3)] there exists a constant $0<\beta<\alpha$ such that
$\|f(t, x_{1},\ldots,x_{r})\|_Y\leq a_1(t)+c_1\sum_{i=1}^{r}\Vert x_{i}\Vert_{X}$
a.e. in $t\in I$, where $a_1\in L^{1/\beta}(I, \mathbb{R}^+)$ and $c_1>0$.
\end{itemize}
	
\item[(H2)] The nonlocal function $h: C(J: X, X)\rightarrow X$ satisfies the following:
\begin{itemize}
\item[(1)] $t\to h(x, y)$ is measurable for all $x,y\in X$;
		
\item[(2)] $\|h(x_{1}, y_{1})-h(x_{2}, y_{2})\|_X\leq k_{2}(t)
\lbrace\|x_{1}-x_{2}\|_X+\|y_{1}-y_{2}\|_X\rbrace$
a.e. on $I$, $k_{2}\in L^{\infty}(I,\mathbb{R}^+)$;
		
\item[(3)] there exists a constant $0<\beta<\alpha$ such that
$\|h(x, y)\|_X\leq a_2(t)+c_2\lbrace\|x\|_X+\|y\|_X\rbrace$
a.e. in $t\in I$ and all $x,y\in X$, where
$a_2\in L^{1/\beta}(\mathbb{R}^+)$ and $c_2>0$.
\end{itemize}
	
\item[(H3)] The multivalued map $U: I\times X
\rightrightarrows P_{f}(T)$ is such that:
\begin{itemize}
\item[(1)] $t\to U(t,x)$ is measurable for all $x\in X$;
		
\item[(2)] $d_H(U(t, x),U(t, y))\leq k_3(t)\|x-y\|_X$ a.e. on $I$,
$k_3\in L^{\infty}(I,\mathbb{R}^+)$;
		
\item[(3)] there exists a constant $0<\beta<\alpha$ such that
$$
\|U(t, x)\|_T=\sup\{\|v\|_T: v\in U(t, x)\}
\leq a_3(t)+c_3\|x\|_X \quad \text{ a.e. in } t\in I,
$$
where $a_3\in L^{1/\beta}(I, \mathbb{R}^+)$ and $c_3>0$.
\end{itemize}
	
\item[(H4)] Functions $g_{i}: I\times X\times T\to \mathbb{R}$,
$i=1,\ldots,r$, are such that:
\begin{itemize}
\item[(1)] the map $t\to g_{i}(t, x, u_{i})$ is measurable
for all $(x, u_{i})\in X\times T$;
	
\item[(2)] $\vert g_{i}(t, x ,u_{i})-g_{i}(t, y ,v_{i})\vert
\leq k^{\prime}_{4}(t)\|x-y\|_X +k^{\prime\prime}_{4}\Vert u_{i}-v_{i}\Vert_{T}$
a.e., $k^{\prime}_{4}\in L^{1}(I,\mathbb{R}^+)$, $k^{\prime\prime}_{4}>0$;
		
\item[(3)] $\vert g_{i}(t, x ,u_{i})\vert\leq a_{4}(t)+b_{4}(t)\Vert x\Vert_{X}
+c_{4}\Vert u_{i}\Vert_{T}$ a.e. $t\in I,
a_4,b_{4}\in L^{1/\beta}(I, \mathbb{R}^+)$, $c_{4}>0$.
\end{itemize}
\end{itemize}

\begin{definition}
\label{Definition 2.4}
A solution of the control system \eqref{eq:1.1}--\eqref{eq:1.3} is defined
to be a vector of functions $(x(\cdot), u_{1}(\cdot),\ldots,u_{r}(\cdot))$
consisting of a trajectory $x\in C(I, X)$ and $r$ multiple controls
$u_{1},\ldots,u_{r}$ $\in L^{1}(I, T)$ satisfying system \eqref{eq:1.1}--\eqref{eq:1.2}
and the inclusion \eqref{eq:1.3} almost everywhere.
\end{definition}

A solution of control system \eqref{eq:1.1}--\eqref{eq:1.2},
\eqref{eq:1.4} can be defined similarly.

\begin{definition}[See \cite{AMA.46,MR3260698,AMA.45}]
\label{Definition 2.5}
A vector of functions $(x,u_{1},\ldots,u_{r})$ is a mild solution of the control
system \eqref{eq:1.1}--\eqref{eq:1.3} iff $x\in C(I,X)$ and there exist
$u_{1},\,\ldots,\,u_{r} \in L^1(I,T)$ such that
$u_{1}(t),\,\ldots,\,u_{r}(t)$ $\in U(t,x(t))$
a.e. in $t\in I$, $x(0)=x_0-h(x(t), B_{r}(t)u_{r}(t))$,
and the following integral equation is satisfied:
\begin{equation*}
x(t)=S_{\alpha}(t)M[x_{0}-h(x(t), B_{r}(t)u_{r}(t))]
+\int_{0}^{t}(t-s)^{\alpha-1}T_{\alpha}(t-s)L^{-1}
f(s, x(s), B_{1}(s)u_{1}(s),\ldots, B_{r-1}(s)u_{r-1}(s))ds,
\end{equation*}
where
$$
S_{\alpha}(t):=\int_{0}^{\infty}M^{-1}\zeta_{\alpha}(\theta)Q(t^{\alpha}\theta)d\theta,
\quad T_{\alpha}(t):=\alpha\int_{0}^{\infty}M^{-1}\theta\zeta_{\alpha}(\theta)
Q(t^{\alpha}\theta)d\theta,
$$
$$
\zeta_{\alpha}(\theta):=\frac{1}{\alpha}\theta^{-1-\frac{1}{\alpha}}
\varpi_{\alpha}(\theta^{-\frac{1}{\alpha}})\geq0,
\quad \varpi_{\alpha}(\theta):=\frac{1}{\pi}\sum_{n=1}^{\infty}(-1)^{n-1}
\theta^{-\alpha n-1}\frac{\Gamma(n\alpha+1)}{n!}\sin (n\pi\alpha),
\ \theta \in ]0, \infty[,
$$
with $\zeta_{\alpha}$ the probability density function defined on $]0, \infty[$,
that is, $\zeta_{\alpha}(\theta)\geq 0$, $\theta\in ]0, \infty[$,
and $\int_{0}^{\infty}\zeta_{\alpha}(\theta)d\theta=1$.
\end{definition}

A similar definition can be introduced for the control
system \eqref{eq:1.1}--\eqref{eq:1.2},\eqref{eq:1.4}.

\begin{remark}[See \cite{AMA.45}]
\label{Remark 2.2}
One has
$\displaystyle \int_0^{\infty}\theta\xi_{\alpha}(\theta)d\theta
=\frac{1}{\Gamma(1+\alpha)}$.
\end{remark}


\subsection{Existence for Multiple Control Systems}

We give existence of solutions for the multiple control systems
\eqref{eq:1.1}--\eqref{eq:1.3} and \eqref{eq:1.1}--\eqref{eq:1.2},\eqref{eq:1.4}.
Let $\Lambda := S(T_\varphi)$. It turns out
that $\Lambda$ is a compact subset of $C(I,X)$ 
and  $\mathcal{T}r_U\subseteq\mathcal{T}
r_{\operatorname{cl} \operatorname{conv} U}\subseteq\Lambda$ \cite{MyID:323}.
Let the set-valued map $\overline{U}: C(I, X) \rightrightarrows 2^{L^{1/\beta}(I, T)}$ be defined by
\begin{equation*}
\overline{U}(x):=\left\{\theta_i: I\to T\text{ measurable}: \theta_i(t)\in U(t, x(t))
\text{ a.e.},~ i=1,\ldots,r\right\},
\quad x\in C(I, X).
\end{equation*}

\begin{theorem}[See \cite{MyID:323}]
\label{Theorem 4.1}
The set $\mathcal{R}_U$ is nonempty and the set
$\mathcal{R}_{\operatorname{cl} \operatorname{conv} U}$ is a compact subset of the space
$C(J,X)\times\omega$-$L^{1/\beta}(I,T)$.
\end{theorem}

\begin{theorem}[See \cite{MyID:323}]
\label{Theorem 5.1}
Let any $(x_*(\cdot),u_{1,*}(\cdot),\ldots,u_{r,*}(\cdot))
\in \mathcal{R}_{\operatorname{cl} \operatorname{conv} U}$.
Then there exists a sequence
$$
(x_n(\cdot), u_{1,n}(\cdot),\ldots,u_{r,n}(\cdot))\in\mathcal{R}_{U},
\quad n\geq1,
$$
such that
\begin{equation}
\label{eq:5.3}
x_n\to x_* ~\text{in}~ C(I, X),
\end{equation}
\begin{equation}
\label{eq:5.4}
u_{i,n}\to u_{i,*}~ \text{in}~ L^{\frac{1}{\beta}}_{\omega}(I, T)~
\text{and}~ \omega\text{-}L^{\frac{1}{\beta}}(I, T),
\end{equation}
\begin{equation}
\label{eq:5.5}
\lim\limits_{n\to \infty}\sup_{0\leq t_1\leq t_2\leq a}
\sum\limits_{i=1}^{r}\biggl\vert\int_{t_1}^{t_2}(g_{i}^{**}(s, x_*(s), u_{i,*}(s))
-g_{i}(s, x_n(s), u_{i,n}(s)))ds
\biggr\vert=0.
\end{equation}
\end{theorem}

\begin{theorem}[See \cite{MyID:323}]
\label{Theorem 5.2}
Problem $(RP)$ has a solution and
\begin{equation}
\label{eq:5.17}
\min\limits_{(x, u_i)\in \mathcal{R}_{\operatorname{cl} \operatorname{conv} U}}J^{**}_{i}(x, u_i)
=\inf\limits_{(x, u_i)\in \mathcal{R}_{U}}J_{i}(x, u_i),
\quad i=1,\ldots,r.
\end{equation}
For any solution $(x_*, u_{1,*},\ldots,u_{r,*})$ of problem $(RP)$,
there exists a minimizing sequence
$$
(x_n, u_{1,n},\ldots,u_{r,n})\in \mathcal{R}_U, \quad n \geq 1,
$$
for problem $(P)$, which converges to $(x_*, u_{1,*},\ldots,u_{r,*})$ in the spaces
$C(I, X)\times \omega\text{-}L^{\frac{1}{\beta}}(I, T)$ and in
$C(I, X)\times L_{\omega}^{\frac{1}{\beta}}(I, T)$, and the following formula holds:
\begin{equation}
\label{eq:5.18}
\lim\limits_{n\to \infty}\sup_{0\leq t_1\leq t_2\leq a}\sum\limits_{i=1}^{r}\biggl\vert
\int_{t_1}^{t_2}(g_{i}^{**}(s, x_*(s), u_{i,*}(s))-g_{i}(s, x_n(s), u_{i,n}(s)))ds
\biggr\vert=0.
\end{equation}
Conversely, if $(x_n, u_{1,n},\ldots,u_{r,n})$, $n\geq1$, is a minimizing sequence for problem $(P)$,
then there is a subsequence $(x_{n_{k}}, u_{1,n_{k}},\ldots,u_{r,n_{k}})$, $k\geq1$, of the
sequence $(x_n, u_{1,n},\ldots,u_{r,n})$, $n\geq1$, and a solution $(x_*, u_{1,*},\ldots,u_{r,*})$
of problem $(RP)$ such that the subsequence $(x_{n_{k}}, u_{1,n_{k}},\ldots,u_{r,n_{k}}),$ $k\geq1$,
converges to $(x_*, u_{1,*},\ldots,u_{r,*})$ in $C(I, X)\times \omega\text{-}L^{\frac{1}{\beta}}(I, T)$
and relation \eqref{eq:5.18} holds for this subsequence $(x_{n_{k}}, u_{1,n_{k}},\ldots,u_{r,n_{k}})$, $k\geq1$.
\end{theorem}

We conclude this survey with the idea that fractional differential equations and fractional optimal 
control are fields of study under strong development. Due to their widespread 
applications  in science and technology, research within the broad area of fractional dynamical 
systems has led to many recent developments that have attracted the attention of a considerable audience 
of scientists. Fractional-order models have the potential to capture nonlocal relations,  
making them more realistic and adequate to describe real-world phenomena. 
In spite of the tremendous number of results in the literature, much remains to be done.


\section*{Acknowledgements}

This work was partially supported by FCT
and CIDMA within project UID/MAT/04106/2013. 
The authors are very grateful to four anonymous referees, 
for their suggestions and invaluable comments.


\bigskip



\end{document}